\documentclass[journal]{IEEEtran}
\usepackage{times}
\usepackage{xcolor}
\usepackage{soul}
\usepackage[utf8]{inputenc}
\usepackage[small]{caption}

\usepackage[cmex10]{amsmath}
\usepackage{epsfig,amsfonts,amssymb,times,subfigure,floatflt,wrapfig,enumerate,tikz,url,bm}
\usepackage{algorithmicx,algpseudocode}
\usepackage{booktabs}
\usepackage{siunitx}
\usepackage{tabu}
\usepackage{graphicx}  
\usepackage{multirow}
\newtheorem{example}{Example}

\newtheorem{definition}{Definition}

\newtheorem{alg}{Algorithm}

\newcommand{\mat}[1]{\bm{#1}}
\newcommand{\ten}[1]{\bm{\mathcal{#1}}}

\ifCLASSOPTIONcompsoc
  \usepackage[nocompress]{cite}
\else
  \usepackage{cite}
\fi

\begin{document}
\title{Fast and Accurate Tensor Completion with Total Variation Regularized Tensor Trains}

\author{Ching-Yun Ko, Kim Batselier,
Wenjian Yu, Senior Member, IEEE,
Ngai Wong, Senior Member, IEEE
\IEEEcompsocitemizethanks{\IEEEcompsocthanksitem Ching-Yun Ko and Ngai Wong are with the Department of Electrical and Electronic Engineering, The University of Hong Kong, Hong Kong. Email: cyko@eee.hku.hk\protect, nwong@eee.hku.hk\protect\\
\IEEEcompsocthanksitem Kim Batselier is with the Delft Center for Systems and Control, Delft University of Technology, Delft, Netherlands. Email: k.batselier@tudelft.nl\protect\\
\IEEEcompsocthanksitem Wenjian Yu is with the Department of Computer Science and Technology, Tsinghua University, Beijing, China. Email: yu-wj@tsinghua.edu.cn\protect\\


}
}

\markboth{Journal of \LaTeX\ Class Files,~Vol.~14, No.~8, August~2015}%
{Shell \MakeLowercase{\textit{et al.}}: Bare Demo of IEEEtran.cls for IEEE Journals}

\maketitle

\IEEEdisplaynontitleabstractindextext

\IEEEpeerreviewmaketitle

\ifCLASSOPTIONcompsoc
\IEEEraisesectionheading{\section{Introduction}\label{sec:introduction}}
\else
\begin{abstract}
We propose a new tensor completion method based on tensor trains. The to-be-completed tensor is modeled as a low-rank tensor train, where we use the known tensor entries and their coordinates to update the tensor train. A novel tensor train initialization procedure is proposed specifically for image and video completion, which is demonstrated to ensure fast convergence of the completion algorithm. The tensor train framework is also shown to easily accommodate Total Variation and Tikhonov regularization due to their low-rank tensor train representations. Image and video inpainting experiments verify the superiority of the proposed scheme in terms of both speed and scalability, where a speedup of up to 155$\times$ is observed compared to state-of-the-art tensor completion methods at a similar accuracy. Moreover, we demonstrate the proposed scheme is especially advantageous over existing algorithms when only tiny portions (say, $1$\%) of the to-be-completed images/videos are known.
\end{abstract}

\begin{IEEEkeywords}
Tensor completion, Tensor-train decomposition, total variation, Image restoration
\end{IEEEkeywords}

%
\IEEEpeerreviewmaketitle

\section{Introduction}
\label{sec:introduction}
\fi

\IEEEPARstart{T}{ensors} are a higher-order generalization of vectors and matrices and have found widespread applications due to their natural representation of real-life multi-way data like images and videos~\cite{xing2012dictionary,lu2016tensor,patwardhan2007video,liu2013tensor,gandy2011tensor}.
Tensor completion generalizes the matrix completion problem, which aims at estimating missing entries from partially revealed data. For example, grayscale images are matrices (two-way tensors) that are indexed by two spatial variables, while color images are essentially three-way tensors with one additional color dimension. Grayscale and color videos are extensions of grayscale and color images by adding one temporal index. Thus image/video recovery tasks are indeed tensor completion problems. Although one can always regard a tensor completion task as multiple matrix completion problems, state-of-the-art matrix completion algorithms such as~\cite{fang2017improved} have rather high computational costs and poor scalability. Moreover, the application of matrix completion methods to tensorial data overlooks a key insight in tensor completion, namely, the low tensor-rank assumption~\cite{liu2013tensor,xu2015parallel,tensorreview,grasedyck2015alternating} inherent to the data. For example, normally every two adjacent frames of a video are shot with a very short time interval, implying that only limited changes are allowed between two adjacent video frames. Similarly, the values among neighboring pixels in an image usually vary slowly. These intuitive low rank ideas have been successfully utilized in research areas such as collaborative filtering~\cite{karatzoglou2010multiverse}, multi-task learning~\cite{romera2013multilinear}, image/video recovery~\cite{xu2015parallel,liu2013tensor,gandy2011tensor} and text analysis~\cite{collins2012tensor,liu2015learning}.  

\subsection{{Related Work}}
\label{subsec:R}
Most existing tensor completion methods are generalizations of matrix completion methods. Traditional matrix completion problems are generally formulated into the construction of a structurally low-rank matrix $\mat{E}$ that has the same observed entries: $$\min_{\mat{E}} {\textrm{rank}(\mat{E})}, \textrm{ \textrm{s.t.} } (\mat{E}-\mat{O})_{\Omega} = \mat{0},$$where $\mat{O}$ represents the observed matrix with zero fillings at the missing entries, and $\Omega$ is the mapping that specifies the locations of known entries. Directly solving the above optimization problem is NP-hard, which resulted in extensive research on solving alternative formulations. Two popular candidates are to minimize the nuclear norm (being the convex envelope of the matrix rank-operator~\cite{candes2009exact,chen2015incoherence}), or to use a factorization method~\cite{wen2012solving} that decomposes the matrix $\mat{E}$ as the product of two small matrices. 
The nuclear norm approach has been generalized to tensor completion problems by unfolding the tensor along its $k$ modes into $k$ matrices and summing over their nuclear norms~\cite{liu2013tensor,collins2012tensor,signoretto2014learning}. Total variation (TV) terms have also been integrated into this method in~\cite{li2017low}. Correspondingly, the factorization method has also been generalized to tensors~\cite{xu2015parallel,tan2014tensor,JiTV}. 

Another way to tackle tensor completion is based on the tensor multi-rank and tubal rank~\cite{kilmer2013third,zhang2014novel,mu2014square}, which are intrinsically defined on three-way tensors such as color images and grayscale videos. It is remarked that multi-rank and tubal rank inspired methods are only applicable when the data are compressible in a t-SVD representation. 
Methods that exploit tensor decomposition formats were also introduced to tensor completion problems in recent years.
In~\cite{suzuki2015convergence} and \cite{zhao2015bayesian}, the authors use the Canonical Polyadic (CP) decomposition for Bayesian tensor estimators. 
Imaizumi et al.~\cite{imaizumi2017tensor} adopted tensor trains (TTs) and solved the optimization problem using a TT Schatten norm~\cite{phien2016efficient} via the alternating direction method of multipliers. 
Bengua et al.~\cite{bengua2017efficient} also combined TT rank optimization with factorization methods introduced in~\cite{xu2015parallel}.
In~\cite{grasedyck2015alternating}, Grasedyck et al. borrowed the traditional tensor completion problem formulation and adopted the tensor train format as the underlying data structure. By updating slices of the tensor cores using parts of the observed entries, the tensor train is completed through alternating least square problems. Wang et al.~\cite{wang2017efficient} designed a tensor completion algorithm by expanding the tensor trains in~\cite{grasedyck2015alternating} to uniform TT-rank tensor rings using random normal distributed values, which yields higher recovery accuracies. However, using tensor rings suffers from two drawbacks: the relatively large storage requirement of a uniform TT-rank structure and a sensitivity of the obtained solution to its random initialization.

\subsection{{Our Contribution}}
\label{subsec:OC}

The motivation for adopting tensor trains in this paper to perform the tensor completion task rather than the CP and Tucker decomposition lies in the fact that determining the CP-rank is NP-hard while the TT-rank is easily determined from an SVD~\cite{oseledets2011tensor}. Also, the Tucker form requires exponentially large storage, which is not as economic as a tensor train. Motivated by these factors, we adopt tensor trains in the tensor completion problem and reformulate the problem as a regression task. 
The unknown completed tensor $\ten{A}$ is thereby interpreted as an underlying regression model, while the observed ``inputs" and ``outputs" of the model are the multi-indices (coordinates) and values of the known tensor entries, respectively. The tensor completion problem is then solved from the following optimization problem
\begin{align*}
\min_{\ten{A}\in\mathcal{S}^{(d)}_{\text{TT}}} \;||\mat{S}^T \;\textrm{vec}(\ten{A}) - \mat{y}||_2^2,&\\
\textrm{\textrm{s.t.}} \;\ten{A}\in\mathcal{S}^{(d)}_{\text{TT}}\text{, and TT-rank}(\ten{A}) = (R_1,&R_2,\ldots,R_{d}),
\end{align*}
where $\mathcal{S}^{(d)}_{\text{TT}}$ denotes the subspace of $d$-order tensor where each element can be represented in the tensor train format with $d$ TT-cores.
The binary matrix $\mat{S}^T$ selects the known entries of $\textrm{vec}(\ten{A})$, where $\textrm{vec}(\ten{A})$ denotes the vectorization of the unknown completed tensor $\ten{A}$. The vector $\mat{y}$ contains the values of the observed tensor entries. Minimizing $||\mat{S}^T \textrm{vec}(\ten{A}) - \mat{y}||_2^2$ therefore enforces that the desired solution needs to have nearly the same observed tensor entries. To regularize the problem, an additional low-rank constraint is added together with the requirement that the desired tensor is represented in the tensor train format, which will be explained in Section~\ref{sec:TNCM}. Moreover, other regularization terms such as Total Variation (TV) and Tikhonov regularization are readily integrated into this tensor train framework. The above problem is solved using an iterative method called the alternating linear scheme that iteratively solves small linear systems. The computationally most expensive steps are a singular value decomposition and QR decomposition. The numerical stability and monotonic convergence of the proposed method is guaranteed by these orthogonal matrix factorizations~\cite{holtz2012alternating,rohwedder2013local}.

The flexibility of our proposed model naturally permits various variants and creates more possibilities for recovery tasks. Particularly, the inputs and outputs in the proposed model can be grouped into batches. This favors parallelization and allows more room in tuning the number of equations and unknowns in least square problems. For example, the updating scheme in~\cite{grasedyck2015alternating,wang2017efficient} can be incorporated as one specific variant of our proposed plain architecture (without TV/ Tikhonov regularizers) by grouping the inputs/outputs into $I_K$ batches when updating the $k$-th core. However, this way of grouping inputs/outputs shows no evidence of balancing the number of equations and unknowns in the resulting least square problems and \cite{grasedyck2015alternating} has demonstrated a consistently inferior performance than~\cite{wang2017efficient}, while we will show in later experiments that the proposed plain TTC model outperforms~\cite{wang2017efficient} in both time and accuracy. Lastly, it is remarked that there is no direct TT format employed in~\cite{bengua2017efficient}. Instead, the TT rank optimization problem formed refers to enforcing the matricizations along $k$ modes have pre-defined rank $=R_k$. This is equivalent to updating the products of the first $k$ tensor train cores and the last $d-k$ cores as two matrices at the same time via a standard factorization method in matrix completions.
The main contributions of this article are:
\begin{itemize}
\item The tensor completion problem is rephrased as a regression task and solved using the alternating linear scheme.
\item Both the selection matrix $\mat{S}$ and tensor $\ten{A}$ are never explicitly formed, which lifts the curse of dimensionality and results in a computationally efficient tensor completion algorithm.
\item A deterministic TT-initialization method for image and video data is provided, which guarantees consistent results and is demonstrated to speed up convergence of our proposed tensor completion algorithm.
\item Both TV and Tikhonov regularization are integrated into the tensor train framework. The addition of these regularization terms comes at almost no additional computational cost. The low TT-rank property of the matrices involved in the TV regularization, especially, ensures a computationally efficient algorithm.
\item We propose a factorization scheme of the physical indices that can be exploited through tensor trains, which results in a better scalability of our method compared to state-of-the-art methods.
\item The proposed method shows a particular effectiveness and efficiency ($155\times$ speedup) in recovering severely damaged (large portions of missing pixels) high-resolution images. A comparable performance to the state-or-the-art methods is demonstrated for small image inpainting tasks.
\end{itemize}
To the best of our knowledge, this is the first time that the Total Variation regularization term is incorporated in tensor train form.
The efficacy of the proposed algorithm is demonstrated by extensive numerical experiments. 
The outline of the article is as follows. Necessary tensor and tensor train preliminaries are briefly introduced in Section~\ref{sec:P}. The proposed tensor train completion model and algorithm are discussed in Section~\ref{sec:TNCM}. Numerical experiments comparing our proposed algorithm with state-of-the-art methods are given in Section~\ref{sec:E}. Finally, conclusions are drawn in Section~\ref{sec:C}.
\section{{Preliminaries}}
\label{sec:P}
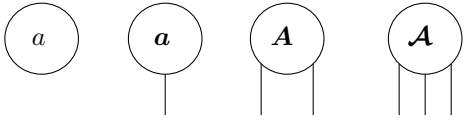
\begin{figure}[tb]
\begin{center}
\ifx\du\undefined
  \newlength{\du}
\fi
\setlength{\du}{4\unitlength}
\begin{tikzpicture}
\pgftransformxscale{1.2000000}
\pgftransformyscale{-1.200000}
\definecolor{dialinecolor}{rgb}{0.000000, 0.000000, 0.000000}
\pgfsetstrokecolor{dialinecolor}
\definecolor{dialinecolor}{rgb}{1.000000, 1.000000, 1.000000}
\pgfsetfillcolor{dialinecolor}
\definecolor{dialinecolor}{rgb}{1.000000, 1.000000, 1.000000}
\pgfsetfillcolor{dialinecolor}
\pgfpathellipse{\pgfpoint{-6.073446\du}{10.779091\du}}{\pgfpoint{2.900000\du}{0\du}}{\pgfpoint{0\du}{2.800000\du}}f
\pgfusepath{fill}
\pgfsetlinewidth{0.100000\du}
\pgfsetdash{}{0pt}
\pgfsetdash{}{0pt}
\definecolor{dialinecolor}{rgb}{0.000000, 0.000000, 0.000000}
\pgfsetstrokecolor{dialinecolor}
\pgfpathellipse{\pgfpoint{-6.073446\du}{10.779091\du}}{\pgfpoint{2.900000\du}{0\du}}{\pgfpoint{0\du}{2.800000\du}}
\pgfusepath{stroke}
\pgfsetlinewidth{0.100000\du}
\pgfsetdash{}{0pt}
\pgfsetdash{}{0pt}
\pgfsetbuttcap
{
\definecolor{dialinecolor}{rgb}{0.000000, 0.000000, 0.000000}
\pgfsetfillcolor{dialinecolor}
\definecolor{dialinecolor}{rgb}{0.000000, 0.000000, 0.000000}
\pgfsetstrokecolor{dialinecolor}
\draw (-4.022836\du,12.758990\du)--(-4.022836\du,16.80\du);
}
\pgfsetlinewidth{0.100000\du}
\pgfsetdash{}{0pt}
\pgfsetdash{}{0pt}
\pgfsetbuttcap
{
\definecolor{dialinecolor}{rgb}{0.000000, 0.000000, 0.000000}
\pgfsetfillcolor{dialinecolor}
\definecolor{dialinecolor}{rgb}{0.000000, 0.000000, 0.000000}
\pgfsetstrokecolor{dialinecolor}
\draw (-8.124055\du,12.758990\du)--(-8.124055\du,16.80\du);
}
\pgfsetlinewidth{0.100000\du}
\pgfsetdash{}{0pt}
\pgfsetdash{}{0pt}
\pgfsetbuttcap
{
\definecolor{dialinecolor}{rgb}{0.000000, 0.000000, 0.000000}
\pgfsetfillcolor{dialinecolor}
\definecolor{dialinecolor}{rgb}{0.000000, 0.000000, 0.000000}
\pgfsetstrokecolor{dialinecolor}
\draw (-6.073446\du,13.579091\du)--(-6.073446\du,16.80\du);
}
\definecolor{dialinecolor}{rgb}{0.000000, 0.000000, 0.000000}
\pgfsetstrokecolor{dialinecolor}
\node[anchor=west] at (-8.20\du,10.6\du){$\ten{A}$};
\definecolor{dialinecolor}{rgb}{1.000000, 1.000000, 1.000000}
\pgfsetfillcolor{dialinecolor}
\pgfpathellipse{\pgfpoint{-36.297515\du}{10.779091\du}}{\pgfpoint{2.900000\du}{0\du}}{\pgfpoint{0\du}{2.800000\du}}
\pgfusepath{fill}
\pgfsetlinewidth{0.100000\du}
\pgfsetdash{}{0pt}
\pgfsetdash{}{0pt}
\definecolor{dialinecolor}{rgb}{0.000000, 0.000000, 0.000000}
\pgfsetstrokecolor{dialinecolor}
\pgfpathellipse{\pgfpoint{-36.297515\du}{10.779091\du}}{\pgfpoint{2.900000\du}{0\du}}{\pgfpoint{0\du}{2.800000\du}}
\pgfusepath{stroke}
\definecolor{dialinecolor}{rgb}{0.000000, 0.000000, 0.000000}
\pgfsetstrokecolor{dialinecolor}
\node[anchor=west] at (-37.80\du,10.80\du){$a$};
\definecolor{dialinecolor}{rgb}{1.000000, 1.000000, 1.000000}
\pgfsetfillcolor{dialinecolor}
\pgfpathellipse{\pgfpoint{-26.560105\du}{10.779091\du}}{\pgfpoint{2.900000\du}{0\du}}{\pgfpoint{0\du}{2.800000\du}}
\pgfusepath{fill}
\pgfsetlinewidth{0.100000\du}
\pgfsetdash{}{0pt}
\pgfsetdash{}{0pt}
\definecolor{dialinecolor}{rgb}{0.000000, 0.000000, 0.000000}
\pgfsetstrokecolor{dialinecolor}
\pgfpathellipse{\pgfpoint{-26.560105\du}{10.779091\du}}{\pgfpoint{2.900000\du}{0\du}}{\pgfpoint{0\du}{2.800000\du}}
\pgfusepath{stroke}
\pgfsetlinewidth{0.100000\du}
\pgfsetdash{}{0pt}
\pgfsetdash{}{0pt}
\pgfsetbuttcap
{
\definecolor{dialinecolor}{rgb}{0.000000, 0.000000, 0.000000}
\pgfsetfillcolor{dialinecolor}
\definecolor{dialinecolor}{rgb}{0.000000, 0.000000, 0.000000}
\pgfsetstrokecolor{dialinecolor}
\draw (-26.560105\du,13.579091\du)--(-26.560105\du,16.80\du);
}
\definecolor{dialinecolor}{rgb}{0.000000, 0.000000, 0.000000}
\pgfsetstrokecolor{dialinecolor}
\node[anchor=west] at (-28.20\du,10.80\du){$\mat{a}$};
\definecolor{dialinecolor}{rgb}{1.000000, 1.000000, 1.000000}
\pgfsetfillcolor{dialinecolor}
\pgfpathellipse{\pgfpoint{-16.948419\du}{10.779091\du}}{\pgfpoint{2.900000\du}{0\du}}{\pgfpoint{0\du}{2.800000\du}}
\pgfusepath{fill}
\pgfsetlinewidth{0.100000\du}
\pgfsetdash{}{0pt}
\pgfsetdash{}{0pt}
\definecolor{dialinecolor}{rgb}{0.000000, 0.000000, 0.000000}
\pgfsetstrokecolor{dialinecolor}
\pgfpathellipse{\pgfpoint{-16.948419\du}{10.779091\du}}{\pgfpoint{2.900000\du}{0\du}}{\pgfpoint{0\du}{2.800000\du}}
\pgfusepath{stroke}
\pgfsetlinewidth{0.100000\du}
\pgfsetdash{}{0pt}
\pgfsetdash{}{0pt}
\pgfsetbuttcap
{
\definecolor{dialinecolor}{rgb}{0.000000, 0.000000, 0.000000}
\pgfsetfillcolor{dialinecolor}
\definecolor{dialinecolor}{rgb}{0.000000, 0.000000, 0.000000}
\pgfsetstrokecolor{dialinecolor}
\draw (-14.897810\du,12.758990\du)--(-14.897810\du,16.80\du);
}
\pgfsetlinewidth{0.100000\du}
\pgfsetdash{}{0pt}
\pgfsetdash{}{0pt}
\pgfsetbuttcap
{
\definecolor{dialinecolor}{rgb}{0.000000, 0.000000, 0.000000}
\pgfsetfillcolor{dialinecolor}
\definecolor{dialinecolor}{rgb}{0.000000, 0.000000, 0.000000}
\pgfsetstrokecolor{dialinecolor}
\draw (-18.999029\du,12.758990\du)--(-18.999029\du,16.80\du);
}
\definecolor{dialinecolor}{rgb}{0.000000, 0.000000, 0.000000}
\pgfsetstrokecolor{dialinecolor}
\node[anchor=west] at (-19.00\du,10.6\du){$\mat{A}$};
\end{tikzpicture}
\caption{Graphical depiction of a scalar $a$, vector $\mat{a}$, matrix $\mat{A}$ and 3-way tensor $\ten{A}$.}
\label{fig:TN}
\end{center}
\end{figure}

Tensors are high-dimensional arrays that generalize vectors and matrices. A $d$-way or $d$-order tensor $\ten{A} \in\mathbb{R}^{I_1\times I_2\times \cdots\times I_d}$ is an array where each entry is indexed by $d$ indices $i_1,i_2,\ldots,i_d$. We use the convention $1\leq i_k\leq I_k$ for $k=1,\ldots,d$. MATLAB notation is used to denote entries of tensors.
In this paper, boldface capital calligraphic letters $\ten{A},\ten{B},\ldots$ are used to denote tensors, boldface capital letters $\mat{A},\mat{B},\ldots$ denote matrices, boldface letters $\mat{a},\mat{b},\ldots$ denote vectors, and Roman letters $a,b,\ldots$ denote scalars. A set of $d$ tensors, like that of a tensor train, is denoted as $\ten{A}^{(1)},\ten{A}^{(2)},\ldots,\ten{A}^{(d)}$.
A useful representation of tensors are tensor network diagrams. These diagrams use a 
graphical depiction of scalars, vectors, matrices, and tensors as introduced in Figure~\ref{fig:TN}, where each node represents a tensor and the number of free edges on each node represents its order. For example, matrices are 2-way tensors, and thus are represented by nodes with two unconnected edges. 
We now give a brief description of some required tensor operations. The generalization of the matrix-matrix multiplication to tensors involves a multiplication of a matrix with a $d$-way tensor along one of its $d$ modes.

\begin{figure*}[t]
\centering
\subfigure{
\begin{minipage}[b]{0.275\linewidth}
\ifx\du\undefined
  \newlength{\du}
\fi
\setlength{\du}{14\unitlength}
\begin{tikzpicture}
\pgftransformxscale{1.00000}
\pgftransformyscale{-1.0000}
\definecolor{dialinecolor}{rgb}{0.000000, 0.000000, 0.000000}
\pgfsetstrokecolor{dialinecolor}
\definecolor{dialinecolor}{rgb}{1.000000, 1.000000, 1.000000}
\pgfsetfillcolor{dialinecolor}
\definecolor{dialinecolor}{rgb}{1.000000, 1.000000, 1.000000}
\pgfsetfillcolor{dialinecolor}
\pgfsetlinewidth{0.050000\du}
\pgfsetdash{}{0pt}
\pgfsetdash{}{0pt}
\pgfsetmiterjoin
\definecolor{dialinecolor}{rgb}{0.000000, 0.000000, 0.000000}
\pgfsetstrokecolor{dialinecolor}
\pgfpathellipse{\pgfpoint{9.788247\du}{5.144124\du}}{\pgfpoint{1.761753\du}{0\du}}{\pgfpoint{0\du}{1.655876\du}}
\pgfusepath{stroke}
\definecolor{dialinecolor}{rgb}{0.000000, 0.000000, 0.000000}
\pgfsetstrokecolor{dialinecolor}
\node at (9.65\du,5.2\du){$\ten{A}$};
\pgfsetlinewidth{0.050000\du}
\pgfsetdash{}{0pt}
\pgfsetdash{}{0pt}
\pgfsetbuttcap
{
\definecolor{dialinecolor}{rgb}{0.000000, 0.000000, 0.000000}
\pgfsetfillcolor{dialinecolor}
\definecolor{dialinecolor}{rgb}{0.000000, 0.000000, 0.000000}
\pgfsetstrokecolor{dialinecolor}
\draw (8.542500\du,6.315005\du)--(8.544052\du,8.5\du);
}
\pgfsetlinewidth{0.050000\du}
\pgfsetdash{}{0pt}
\pgfsetdash{}{0pt}
\pgfsetbuttcap
{
\definecolor{dialinecolor}{rgb}{0.000000, 0.000000, 0.000000}
\pgfsetfillcolor{dialinecolor}
\definecolor{dialinecolor}{rgb}{0.000000, 0.000000, 0.000000}
\pgfsetstrokecolor{dialinecolor}
\draw (11.033995\du,6.315005\du)--(11.029923\du,8.5\du);
}
\pgfsetlinewidth{0.050000\du}
\pgfsetdash{}{0pt}
\pgfsetdash{}{0pt}
\pgfsetbuttcap
{
\definecolor{dialinecolor}{rgb}{0.000000, 0.000000, 0.000000}
\pgfsetfillcolor{dialinecolor}
\definecolor{dialinecolor}{rgb}{0.000000, 0.000000, 0.000000}
\pgfsetstrokecolor{dialinecolor}
\draw (9.317034\du,6.703578\du)--(9.317034\du,8.5\du);
}
\pgfsetlinewidth{0.050000\du}
\pgfsetdash{}{0pt}
\pgfsetdash{}{0pt}
\pgfsetbuttcap
{
\definecolor{dialinecolor}{rgb}{0.000000, 0.000000, 0.000000}
\pgfsetfillcolor{dialinecolor}
\definecolor{dialinecolor}{rgb}{0.000000, 0.000000, 0.000000}
\pgfsetstrokecolor{dialinecolor}
\draw (10.211006\du,6.770201\du)--(10.211006\du,8.5\du);
}
\definecolor{dialinecolor}{rgb}{0.000000, 0.000000, 0.000000}
\pgfsetstrokecolor{dialinecolor}
\node[anchor=west] at (8.1\du,9.1\du){$I_1$};
\definecolor{dialinecolor}{rgb}{0.000000, 0.000000, 0.000000}
\pgfsetstrokecolor{dialinecolor}
\node[anchor=west] at (8.9\du,9.1\du){$I_2$};
\definecolor{dialinecolor}{rgb}{0.000000, 0.000000, 0.000000}
\pgfsetstrokecolor{dialinecolor}
\node[anchor=west] at (9.8\du,9.1\du){$I_3$};
\definecolor{dialinecolor}{rgb}{0.000000, 0.000000, 0.000000}
\pgfsetstrokecolor{dialinecolor}
\node[anchor=west] at (10.7\du,9.1\du){$I_4$};
\pgfsetdash{}{0pt}
\pgfsetdash{}{0pt}
\pgfsetbuttcap
\pgfsetmiterjoin
\pgfsetlinewidth{0.050000\du}
\pgfsetbuttcap
\pgfsetmiterjoin
\pgfsetdash{}{0pt}
\definecolor{dialinecolor}{rgb}{1.000000, 1.000000, 1.000000}
\pgfsetfillcolor{dialinecolor}
\fill (15.292782\du,6.257627\du)--(16.898207\du,6.257627\du)--(16.898207\du,5.744783\du)--(18.503631\du,6.770471\du)--(16.898207\du,7.796159\du)--(16.898207\du,7.283315\du)--(15.292782\du,7.283315\du)--cycle;
\definecolor{dialinecolor}{rgb}{0.000000, 0.000000, 0.000000}
\pgfsetstrokecolor{dialinecolor}
\draw (14.0\du,6.257627\du)--(16.0\du,6.257627\du)--(16.0\du,5.744783\du)--(17.503631\du,6.770471\du)--(16.0\du,7.796159\du)--(16.0\du,7.283315\du)--(14.0\du,7.283315\du)--cycle;
\end{tikzpicture}
\end{minipage}}
\subfigure{
\begin{minipage}[b]{0.4\linewidth}
\ifx\du\undefined
  \newlength{\du}
\fi
\setlength{\du}{4.5\unitlength}
\begin{tikzpicture}
\pgftransformxscale{1.0000}
\pgftransformyscale{-1.0000}
\definecolor{dialinecolor}{rgb}{0.000000, 0.000000, 0.000000}
\pgfsetstrokecolor{dialinecolor}
\definecolor{dialinecolor}{rgb}{1.000000, 1.000000, 1.000000}
\pgfsetfillcolor{dialinecolor}
\definecolor{dialinecolor}{rgb}{1.000000, 1.000000, 1.000000}
\pgfusepath{fill}
\pgfsetlinewidth{0.100000\du}
\pgfsetdash{}{0pt}
\pgfsetdash{}{0pt}
\definecolor{dialinecolor}{rgb}{0.000000, 0.000000, 0.000000}
\pgfsetstrokecolor{dialinecolor}
\pgfpathellipse{\pgfpoint{-4.7267\du}{9.997200\du}}{\pgfpoint{3.550000\du}{0\du}}{\pgfpoint{0\du}{3.450000\du}}
\pgfusepath{stroke}
\pgfsetlinewidth{0.100000\du}
\pgfsetdash{}{0pt}
\pgfsetdash{}{0pt}
\pgfsetbuttcap
{
\definecolor{dialinecolor}{rgb}{0.000000, 0.000000, 0.000000}
\pgfsetfillcolor{dialinecolor}
\definecolor{dialinecolor}{rgb}{0.000000, 0.000000, 0.000000}
\pgfsetstrokecolor{dialinecolor}
\draw (-1.1767\du,10.036600\du)--(3.196473\du,10.036600\du);
}
\definecolor{dialinecolor}{rgb}{0.000000, 0.000000, 0.000000}
\pgfsetstrokecolor{dialinecolor}
\node[anchor=west] at (-6.3\du,20\du){$I_1$};
\definecolor{dialinecolor}{rgb}{1.000000, 1.000000, 1.000000}
\pgfsetfillcolor{dialinecolor}
\pgfpathellipse{\pgfpoint{6.746473\du}{9.997200\du}}{\pgfpoint{3.550000\du}{0\du}}{\pgfpoint{0\du}{3.450000\du}}
\pgfusepath{fill}
\pgfsetlinewidth{0.100000\du}
\pgfsetdash{}{0pt}
\pgfsetdash{}{0pt}
\definecolor{dialinecolor}{rgb}{0.000000, 0.000000, 0.000000}
\pgfsetstrokecolor{dialinecolor}
\pgfpathellipse{\pgfpoint{6.746473\du}{9.997200\du}}{\pgfpoint{3.550000\du}{0\du}}{\pgfpoint{0\du}{3.450000\du}}
\pgfusepath{stroke}
\pgfsetlinewidth{0.100000\du}
\pgfsetdash{}{0pt}
\pgfsetdash{}{0pt}
\pgfsetbuttcap
{
\definecolor{dialinecolor}{rgb}{0.000000, 0.000000, 0.000000}
\pgfsetfillcolor{dialinecolor}
\definecolor{dialinecolor}{rgb}{0.000000, 0.000000, 0.000000}
\pgfsetstrokecolor{dialinecolor}
\draw (6.746473\du,13.447200\du)--(6.691773\du,17.927200\du);
}
\definecolor{dialinecolor}{rgb}{0.000000, 0.000000, 0.000000}
\pgfsetstrokecolor{dialinecolor}
\node[anchor=west] at (5.\du,20\du){$I_2$};
\definecolor{dialinecolor}{rgb}{1.000000, 1.000000, 1.000000}
\pgfsetfillcolor{dialinecolor}
\pgfpathellipse{\pgfpoint{18.147073\du}{9.997200\du}}{\pgfpoint{3.550000\du}{0\du}}{\pgfpoint{0\du}{3.450000\du}}
\pgfusepath{fill}
\pgfsetlinewidth{0.100000\du}
\pgfsetdash{}{0pt}
\pgfsetdash{}{0pt}
\definecolor{dialinecolor}{rgb}{0.000000, 0.000000, 0.000000}
\pgfsetstrokecolor{dialinecolor}
\pgfpathellipse{\pgfpoint{18.147073\du}{9.997200\du}}{\pgfpoint{3.550000\du}{0\du}}{\pgfpoint{0\du}{3.450000\du}}
\pgfusepath{stroke}
\pgfsetlinewidth{0.100000\du}
\pgfsetdash{}{0pt}
\pgfsetdash{}{0pt}
\pgfsetbuttcap
{
\definecolor{dialinecolor}{rgb}{0.000000, 0.000000, 0.000000}
\pgfsetfillcolor{dialinecolor}
\definecolor{dialinecolor}{rgb}{0.000000, 0.000000, 0.000000}
\pgfsetstrokecolor{dialinecolor}
\draw (21.697073\du,9.997200\du)--(26.070273\du,9.997200\du);
}
\pgfsetlinewidth{0.100000\du}
\pgfsetdash{}{0pt}
\pgfsetdash{}{0pt}
\pgfsetbuttcap
{
\definecolor{dialinecolor}{rgb}{0.000000, 0.000000, 0.000000}
\pgfsetfillcolor{dialinecolor}
\definecolor{dialinecolor}{rgb}{0.000000, 0.000000, 0.000000}
\pgfsetstrokecolor{dialinecolor}
\draw (10.296473\du,10.036600\du)--(14.597073\du,10.036600\du);
}
\definecolor{dialinecolor}{rgb}{0.000000, 0.000000, 0.000000}
\pgfsetstrokecolor{dialinecolor}
\node[anchor=west] at (16.5\du,20\du){$I_3$};
\definecolor{dialinecolor}{rgb}{1.000000, 1.000000, 1.000000}
\pgfsetfillcolor{dialinecolor}
\pgfpathellipse{\pgfpoint{29.620273\du}{9.997200\du}}{\pgfpoint{3.550000\du}{0\du}}{\pgfpoint{0\du}{3.450000\du}}
\pgfusepath{fill}
\pgfsetlinewidth{0.100000\du}
\pgfsetdash{}{0pt}
\pgfsetdash{}{0pt}
\definecolor{dialinecolor}{rgb}{0.000000, 0.000000, 0.000000}
\pgfsetstrokecolor{dialinecolor}
\pgfpathellipse{\pgfpoint{29.620273\du}{9.997200\du}}{\pgfpoint{3.550000\du}{0\du}}{\pgfpoint{0\du}{3.450000\du}}
\pgfusepath{stroke}
\definecolor{dialinecolor}{rgb}{0.000000, 0.000000, 0.000000}
\pgfsetstrokecolor{dialinecolor}
\node[anchor=west] at (28\du,20\du){$I_4$};
\definecolor{dialinecolor}{rgb}{0.000000, 0.000000, 0.000000}
\pgfsetstrokecolor{dialinecolor}
\node[anchor=west] at (-8\du,9.554100\du){$\ten{A}^{(1)}$};
\definecolor{dialinecolor}{rgb}{0.000000, 0.000000, 0.000000}
\pgfsetstrokecolor{dialinecolor}
\node[anchor=west] at (4\du,9.554100\du){$\ten{A}^{(2)}$};
\definecolor{dialinecolor}{rgb}{0.000000, 0.000000, 0.000000}
\pgfsetstrokecolor{dialinecolor}
\node[anchor=west] at (15.\du,9.554100\du){$\ten{A}^{(3)}$};
\definecolor{dialinecolor}{rgb}{0.000000, 0.000000, 0.000000}
\pgfsetstrokecolor{dialinecolor}
\node[anchor=west] at (26.5\du,9.554100\du){$\ten{A}^{(4)}$};
\definecolor{dialinecolor}{rgb}{0.000000, 0.000000, 0.000000}
\pgfsetstrokecolor{dialinecolor}
\node[anchor=west] at (11.697173\du,3.\du){$R_1$};
\definecolor{dialinecolor}{rgb}{0.000000, 0.000000, 0.000000}
\pgfsetstrokecolor{dialinecolor}
\node[anchor=west] at (-1\du,8\du){$R_2$};
\definecolor{dialinecolor}{rgb}{0.000000, 0.000000, 0.000000}
\pgfsetstrokecolor{dialinecolor}
\node[anchor=west] at (10\du,8\du){$R_3$};
\definecolor{dialinecolor}{rgb}{0.000000, 0.000000, 0.000000}
\pgfsetstrokecolor{dialinecolor}
\node[anchor=west] at (21\du,8\du){$R_4$};
\pgfsetlinewidth{0.100000\du}
\pgfsetdash{}{0pt}
\pgfsetdash{}{0pt}
\pgfsetmiterjoin
\pgfsetbuttcap
{
\definecolor{dialinecolor}{rgb}{0.000000, 0.000000, 0.000000}
\pgfsetfillcolor{dialinecolor}
\definecolor{dialinecolor}{rgb}{0.000000, 0.000000, 0.000000}
\pgfsetstrokecolor{dialinecolor}
\pgfpathmoveto{\pgfpoint{-8.2767\du}{9.997200\du}}\,
\pgfpathcurveto{\pgfpoint{-19.669327\du}{2.594600\du}}{\pgfpoint{46.340273\du}{2.983600\du}}{\pgfpoint{33.170273\du}{9.997200\du}}
\pgfusepath{stroke}
}
\pgfsetlinewidth{0.100000\du}
\pgfsetdash{}{0pt}
\pgfsetdash{}{0pt}
\pgfsetbuttcap
{
\definecolor{dialinecolor}{rgb}{0.000000, 0.000000, 0.000000}
\pgfsetfillcolor{dialinecolor}
\definecolor{dialinecolor}{rgb}{0.000000, 0.000000, 0.000000}
\pgfsetstrokecolor{dialinecolor}
\draw (-4.7267\du,13.447200\du)--(-4.7544\du,17.890607\du);
}
\pgfsetlinewidth{0.100000\du}
\pgfsetdash{}{0pt}
\pgfsetdash{}{0pt}
\pgfsetbuttcap
{
\definecolor{dialinecolor}{rgb}{0.000000, 0.000000, 0.000000}
\pgfsetfillcolor{dialinecolor}
\definecolor{dialinecolor}{rgb}{0.000000, 0.000000, 0.000000}
\pgfsetstrokecolor{dialinecolor}
\draw (18.147073\du,13.447200\du)--(18.150000\du,17.900000\du);
}
\pgfsetlinewidth{0.100000\du}
\pgfsetdash{}{0pt}
\pgfsetdash{}{0pt}
\pgfsetbuttcap
{
\definecolor{dialinecolor}{rgb}{0.000000, 0.000000, 0.000000}
\pgfsetfillcolor{dialinecolor}
\definecolor{dialinecolor}{rgb}{0.000000, 0.000000, 0.000000}
\pgfsetstrokecolor{dialinecolor}
\draw (29.620273\du,13.447200\du)--(29.535607\du,17.860607\du);
}
\end{tikzpicture}
\end{minipage}}
\caption{Graphical depiction of the tensor train decomposition of a 4-way tensor $\ten{A}$.}
\label{fig:tensor2mps}
\end{figure*}
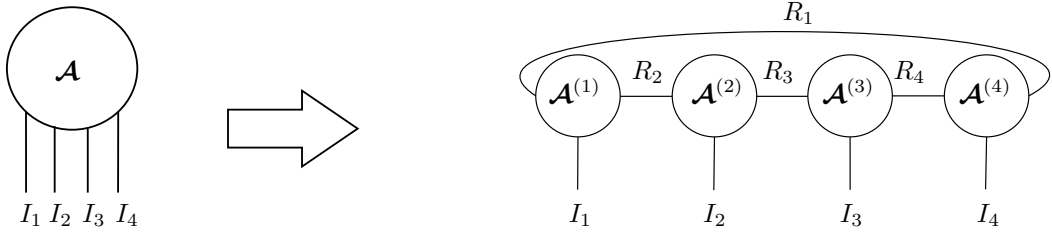

\begin{definition}(\cite[p.~460]{tensorreview})
The $k$-mode product of a tensor~\mbox{$\ten{A}\in\mathbb{R}^{I_1\times\cdots \times I_d}$} with a matrix $\bm{U}\in\mathbb{R}^{J \times I_k}$ is denoted \mbox{$\ten{B}=\ten{A}\, {\times_k}\, \bm{U}$} and defined by%
\begin{align*}
\nonumber \ten{B}(i_1,\cdots,i_{k-1},j,i_{k+1}, \cdots, i_d) &= \hfill \\
\sum\limits_{i_k=1}^{I_k}  \mat{U}(j,i_k) \ten{A}(i_1,\cdots,i_{k-1},i_k,&i_{k+1},\cdots,i_d),
\end{align*}
where $\ten{B}\in\mathbb{R}^{I_1\times\cdots \times I_{k-1}\times J \times I_{k+1}\times\cdots\times I_d}$.
\end{definition}
The proposed method also requires the knowledge of the matrix Khatri-Rao product, which is crucial in the construction of the input matrix $\mat{S}$. 
\begin{definition}
If $\mat{A} \in \mathbb{R}^{N_1 \times M}$ and~$\mat{C} \in \mathbb{R}^{N_2 \times M}$, then their Khatri-Rao product $\mat{A} \odot \mat{C}$ is the $N_1N_2 \times M$ matrix
\begin{align*}
\begin{pmatrix}
\mat{A}(:,1)\otimes \mat{C}(:,1) & \cdots & \mat{A}(:,M)\otimes \mat{C}(:,M) \\
\end{pmatrix},
\end{align*}
where $\otimes$ denotes the standard Kronecker product.
\end{definition}%
Two other basic operations include tensor reshaping and vectorization.
\begin{definition} We adopt the MATLAB reshape operator ``reshape($\ten{A},[n_1,n_2,n_3 \cdots])$", which reshapes the $d$-way tensor $\ten{A}$ with column-wise ordering preserved into a tensor with dimensions $n_1 \times n_2 \times \cdots \times n_d$. The total number of elements of $\ten{A}$ must be the same as~$n_1\times n_2 \times \cdots \times n_d$.
\end{definition}
\begin{definition}
The vectorization $\textrm{vec}(\ten{A})$ of a tensor $\ten{A}$ is the vector obtained from concatenating all tensor entries into one column vector.
\end{definition}
More details on these tensor operations can be found in~\cite[p.~459]{tensorreview}. The mapping between the index $i$ of the vector $\textrm{vec}(\ten{A})$ and the multi-index $[i_1,i_2,\ldots,i_d]$ of the corresponding tensor $\ten{A}$ is bijective. Herein we adopt the mapping convention%
\begin{align}
i &= 
i_1+\sum_{k=2}^d{(i_{k}-1)\prod_{p=1}^{k-1}I_{p}}
\label{eqn:single2multi}
\end{align}
to convert the multi-index $[i_1,i_2,\ldots,i_d]$ into the linear index $i$. Before going into the formal definition of a tensor train decomposition, we give an intuitive interpretation of this process in Figure~\ref{fig:tensor2mps}. For a given 4-way tensor as the one on the left-hand side, we analogize it to a large train cabin with four wheels (legs), while a tensor train decomposition separates the train cabin to four small connected cabins as shown on the right-hand side. 
Using a tensor train format can reduce the storage complexity of a tensor significantly. The storage of a tensor train of a $d$-way tensor $\ten{Z}\in\mathbb{R}^{I\times I\times\cdots\times I}$ requires approximately $dIR^2$ elements, where $R$ is the maximum TT-rank, compared to the conventional $I^d$ storage requirement.
The formal definition is given as follows:

\begin{definition}
A tensor train of a vector $\mat{a} \in \mathbb{R}^{I_1I_2\cdots I_d}$ is defined as the set of $d$ 3-way tensors $\ten{A}^{(1)},\ldots,\ten{A}^{(d)}$ such that $\mat{a}(i)$ equals
{
\small
\begin{align}
\label{eqn:mpodef}
 \sum_{r_1=1}^{R_1}\cdots\sum_{r_d=1}^{R_d} \ten{A}^{(1)}(r_1,i_1,r_2)\,\ten{A}^{(2)}(r_2,i_2,r_3)\cdots \ten{A}^{(d)}(r_d,i_d,r_{1}),%
\end{align}}%
where $i$ is related to $[i_1,i_2,\ldots,i_d]$ via \eqref{eqn:single2multi} and $r_1,r_2,\ldots,r_{d}$ are auxiliary indices that are summed over. The dimensions $R_1,R_2,\ldots,R_{d}$ of these auxiliary indices are called the tensor train ranks (TT-ranks).
\end{definition}
The summations over the auxiliary indices are represented in Figure~\ref{fig:tensor2mps} by the connected edges between the different nodes.
The second auxiliary index of $\ten{A}^{(d)}$ is $R_1$, which ensures that the summation in \eqref{eqn:mpodef} results in a scalar. When $R_1 >1$, the tensor train is also called a tensor ring or matrix product state with periodic boundary conditions~\cite{zhao2018learning,Orus2013}. Throughout this article we always choose \mbox{$R_1=1$}.


The notion of a site-$k$-mixed canonical tensor train is very useful when implementing our proposed algorithm. In order to be able to define this notion, we first need to introduce left and right-orthogonal tensor train cores.
\begin{definition}(\cite[p.~A689]{holtz2012alternating})
A tensor train core $\ten{A}^{(k)}$ is left-orthogonal if it can be reshaped into an $R_{k}I_k\times R_{k+1}$ matrix $\mat{A}$ for which
\begin{align*}
\mat{A}^T\mat{A} &=\mat{I}_{R_{k+1}}
\end{align*}
applies. Similarly, a tensor train core $\ten{A}^{(k)}$ is right-orthogonal if it can be reshaped into an $R_{k}\times I_kR_{k+1}$ matrix $\mat{A}$ for which
\begin{align*}
\mat{A}\mat{A}^T &=\mat{I}_{R_k}
\end{align*}
applies.
\end{definition}
\begin{definition}(Site-$k$-mixed-canonical tensor train)
A tensor train is in site-$k$-mixed-canonical form~\cite{SCHOLLWOCK201196} when all TT-cores $\ten{A}^{(l)} \,(1\leq l \leq k-1)$ are left-orthogonal and TT-cores $\ten{A}^{(l)}\,(k+1 \leq l \leq d)$ are right-orthogonal.
\end{definition}
One advantage of a tensor $\ten{A}$ being in a site-$k$-mixed-canonical form is that its Frobenius norm is easily computed as
\begin{align}
\label{eqn:norm}
{||\ten{A}||}_F^2&={||\ten{A}^{(k)}||}_F^2 = \textrm{vec}(\ten{A}^{(k)})^T \textrm{vec}(\ten{A}^{(k)}).
\end{align}

\section{Methodology}
\label{sec:TNCM}
We first explain the basic idea of our proposed method without any TV or Tikhonov regularization in Sections~\ref{subsec:basicidea} up to \ref{subsec:ALS}. Discussions on how the tensor dimensions can be factorized and the TT-ranks chosen are given in Section~\ref{subsec:factordimensions&TTranks}. The inclusion of both TV and Tikhonov regularization are discussed in Sections~\ref{subsec:TV} and \ref{subsec:Tikhonov}, respectively.
\subsection{Basic Idea}
\label{subsec:basicidea}
The proposed tensor completion method intrinsically relies on solving an underdetermined linear system under a low TT-rank constraint. For a given set of $N$ multi-indices $[i_1,i_2,\ldots,i_d]$ and a corresponding vector of observed tensor entries $\mat{y} \in \mathbb{R}^N$, the goal is to obtain a tensor $\ten{A} \in \mathbb{R}^{I_1 \times \cdots \times I_d}$ that contains the same tensor entries. Equivalently, we form the optimization problem
\begin{align}
\label{eqn:ude}
\min_{\ten{A}\in\mathcal{S}^{(d)}_{\text{TT}}} \;||\mat{S}^T \;\textrm{vec}(\ten{A}) - \mat{y}||_2^2,&\\
\nonumber\textrm{\textrm{s.t.}} \;\ten{A}\in\mathcal{S}^{(d)}_{\text{TT}}\text{, and TT-rank}(\ten{A}) = (R_1,&R_2,\ldots,R_{d}),
\end{align}
where $\mat{S}^T \in\mathbb{R}^{N\times I_1I_2I_3\cdots I_d}$ is the row selection matrix corresponding with the known multi-indices. In other words, we want $\mat{S}^T\textrm{vec}(\ten{A})\in\mathbb{R}^{N\times 1}$ to be as close as possible to the observed entries $\mat{y}$, under the constraint that $\textrm{vec}(\ten{A})$ has a low-rank tensor train representation. We stress that the tensor $\ten{A}$ is never computed explicitly in the proposed algorithm but is stored in its tensor train format instead. This is possible by exploiting the fact that $\mat{S}$ can be written as the Khatri-Rao product of $d$ smaller binary matrices 
\begin{align}
\label{eqn:khatriS}
\mat{S} &= \mat{S}^{(d)}\odot\mat{S}^{(d-1)}\odot\ldots\odot\mat{S}^{(1)}.
\end{align}
Equation \eqref{eqn:khatriS} also implies that the matrix $\mat{S}$ never needs to be explicitly kept in memory. By storing the factor matrices $\mat{S}^{(1)},\ldots,\mat{S}^{(d)}$ instead, the storage cost for $\mat{S}$ is reduced from $N(I_1\cdots I_d)$ down to $N(I_1+\cdots+I_d)$. The decomposition of $\mat{S}$ follows from the following definition.




\begin{definition}
\label{def:selvec}
For a tensor entry $\ten{A}(i_1,i_2,\ldots,i_d)$, the corresponding selection vector $\mat{s}_{[i_1,i_2,\ldots,i_d]}$ is defined as
\begin{align}
\label{eqn:selvec}
\mat{s}_{[i_1,i_2,\ldots,i_d]}:=\mat{e_{i_d}}^{(d)}\otimes \cdots \otimes \mat{e_{i_2}}^{(2)}\otimes \mat{e_{i_1}}^{(1)},
\end{align}
where $\mat{e_{i_k}}^{(k)}\in\mathbb{R}^{I_k}\; (k=1,2,\ldots,d)$ is the $i_k$-th standard basis vector.
\end{definition}
One can verify that 
$$\mat{s}_{[i_1,i_2,\ldots,i_d]}^T\,\textrm{vec}(\ten{A})=\ten{A}(i_1,i_2,\ldots,i_d).$$
Note that the order of the Kronecker products is reversed to be consistent with the index mapping~\eqref{eqn:single2multi}. We now define the $n$-th column of the $N\times I_k$ selection matrix $\mat{S}^{(k)}$ as the standard basis vector $\mat{e_{i_k}}^{(k)}$ of the $n$-th observed entry. Equation~\eqref{eqn:khatriS} then follows from the concatenation of \eqref{eqn:selvec} for each known multi-index.
\begin{example}
\label{ex:1}
We use a small example to illustrate Definition~\ref{def:selvec}. Consider a $3\times 4\times 2$ tensor
$\ten{A}$ with only 3 of the entries observed
. The multi-indices of the three observed entries 
are
\begin{align*}
[2,1,2],[1,3,1],[3,4,2].
\end{align*}
The corresponding selection vectors are then given by
\begin{align*}
\mat{s}_{[2,1,2]}&=\mat{e_2}^{(3)}\otimes \mat{e_1}^{(2)}\otimes \mat{e_2}^{(1)},\\
\mat{s}_{[1,3,1]}&=\mat{e_1}^{(3)}\otimes \mat{e_3}^{(2)}\otimes \mat{e_1}^{(1)},\\
\mat{s}_{[3,4,2]}&=\mat{e_2}^{(3)}\otimes \mat{e_4}^{(2)}\otimes \mat{e_3}^{(1)}.
\end{align*}
The corresponding selection matrices for each mode are then
\begin{align*}
\mat{S}^{(1)} &= 
\begin{pmatrix}
\mat{e_2}^{(1)} & \mat{e_1}^{(1)} & \mat{e_3}^{(1)}
\end{pmatrix}  \in\mathbb{R}^{3\times 3},\\
\mat{S}^{(2)} &= 
\begin{pmatrix}
\mat{e_1}^{(2)} & \mat{e_3}^{(2)} & \mat{e_4}^{(2)}
\end{pmatrix} \in\mathbb{R}^{4\times 3},\\
\mat{S}^{(3)} &= 
\begin{pmatrix}
\mat{e_2}^{(3)} & \mat{e_1}^{(3)} & \mat{e_2}^{3)}
\end{pmatrix} \in\mathbb{R}^{2\times 3}.
\end{align*}
\end{example}
In our regression task, the $d$ matrices $\mat{S}^{(1)},\mat{S}^{(2)},\ldots,\mat{S}^{(d)}$ from~\eqref{eqn:khatriS} act as the inputs to the unknown model and the output is the vector of $N$ observed tensor entries
\begin{align*}
\mat{y} &= (y_1,y_2,\ldots,y_N)^T,
\end{align*}
where typically $N \ll I_1I_2\cdots I_d$. 
\subsection{{Tensor Train Initialization}}
\label{subsec:Init}
The proposed tensor completion algorithm solves~\eqref{eqn:ude} iteratively using an alternating linear scheme (ALS). Starting from an initial guess for the tensor train of $\textrm{vec}(\ten{A})$, the ALS updates each tensor train core for a predefined number of iterations or until convergence. Each tensor train core update is achieved by solving a relatively small (compared to the original tensor size) least squares problem. A good initial guess is therefore of crucial importance to speed up convergence. Through extensive tests, we observed that the following heuristical initialization method for images and videos is effective in terms of convergence speed. 

Suppose $\ten{V}\in\mathbb{R}^{I_1\times I_2\times\cdots\times I_d}$ is a tensor with missing entries. Before converting $\ten{V}$ into a tensor train, the goal is to fill up the missing entries using information from the observed entries through two interpolation steps. First, each dimension of the tensor $\ten{V}$ is resized by a factor $h$ using a box-shaped interpolation kernel. The resulting tensor is denoted \mbox{$\ten{W}\in\mathbb{R}^{\lfloor\frac{I_1}{h}\rfloor\times \lfloor\frac{I_2}{h}\rfloor\times\cdots\times \lfloor\frac{I_d}{h}\rfloor}$} and its entries are then used to construct a tensor of the original size $\ten{X}\in\mathbb{R}^{I_1\times I_2\times\cdots\times I_d}$ through cubic kernel interpolation. Alternatively, one can use max or average pooling together with interpolation to achieve a similar effect. Note that for colour images and videos the colour dimension is not resized during the whole initialization procedure. Finally, a tensor train with given TT-ranks is computed from $\ten{X}$ by a modified version of the TT-SVD algorithm~\cite[p.~2301]{oseledets2011tensor}. The TT-SVD algorithm decomposes a tensor to its tensor train format by consecutive reshapings and singular value decomposition (SVD). The modification is made in line $5$ of the TT-SVD algorithm, where instead of using the original truncation parameter $\delta$, each SVD is truncated to the prescribed TT-rank $R_2,\ldots,R_{d}$. Alternatively, one can use a Krylov subspace method (`svds' command in MATLAB) to determine the desired truncated SVD. Using the TT-SVD algorithm to obtain the initial estimate of tensor train also implies that the tensor train will be in site-$d$-mixed-canonical form.

\subsection{{Alternating Linear Scheme}}

We now derive the least squares problem for updating each tensor train core during the ALS. The main motivation for solving~\eqref{eqn:ude} in tensor train form is the reduction in computational cost. Indeed, we will show how updating each tensor train core $\ten{A}^{(k)}$ during the ALS has a computational cost of $O(N (R_kI_kR_{k+1})^2)$ flops, whereas in vector form the computational cost would take approximately $O(N^2I^d)$ flops. In addition, by specifying small TT-ranks one effectively regularizes the problem since the underdetermined system $\mat{S}^T\textrm{vec}(\ten{A})=\mat{y}$ will typically have an infinite number of solutions. In what follows, $\mat{s}^{(k)}_l\in\mathbb{R}^{I_k\times 1} (1\leq l \leq N)$ denotes the $l$-th column of the matrix $\mat{S}^{(k)}$. We further define the following useful auxiliary notations
\begin{align*}
\mat{a}_{<k,l}^T &:= (\ten{A}^{(1)} \times_2 \mat{s}^{(1)T}_l)  \ldots (\ten{A}^{(k-1)} \times_2 \mat{s}^{(k-1)T}_l) \in\mathbb{R}^{R_k},\\
\mat{a}_{>k,l} &:= (\ten{A}^{(k+1)} \times_2 \mat{s}^{(k+1)T}_l ) \ldots (\ten{A}^{(d)} \times_2 \mat{s}^{(d)T}_l) \in\mathbb{R}^{R_{k+1}},
\end{align*}
for \mbox{$k=2,\ldots,d-1$}. Per definition $\mat{a}_{<1,l}=\mat{a}_{>d,l}=1$. The $l$-th observed entry $\mat{y}(l)$ can then be written as 
\begin{align}  
\label{eq:ybba}
\mat{y}(l) &=(\mat{a}_{>k,l}^T\otimes\mat{s}^{(k)T}_l\otimes\mat{a}_{<k,l})\, \textrm{vec}(\ten{A}^{(k)}).
\end{align}
The proof of equation~\eqref{eq:ybba} resembles that in~\cite[Theorem~4.1]{Batselier2017TensorNA}. Writing out~\eqref{eq:ybba} for all $N$ observed entries results in the following linear system
\begin{align}
\label{eq:le}
\mat{y}
&=
\begin{pmatrix}
\mat{a}_{>k,1}^T\otimes\mat{s}^{(k)T}_1\otimes\mat{a}_{<k,1}\\
\mat{a}_{>k,2}^T\otimes\mat{s}^{(k)T}_2\otimes\mat{a}_{<k,2}\\
\vdots\\
\mat{a}_{>k,N}^T\otimes\mat{s}^{(k)T}_N\otimes\mat{a}_{<k,N}
\end{pmatrix}
\textrm{vec}(\ten{A}^{(k)}),
\end{align}
where the matrix is $N \times R_kI_kR_{k+1}$. 
Solving \eqref{eq:le} requires $O(N (R_kI_kR_{k+1})^2)$ flops. 
If the TT-ranks $R_k,R_{k+1}$ are chosen such that $N \geq R_kI_kR_{k+1}$, one can solve the normal equations of~\eqref{eq:le} instead with a computational complexity of $O((R_kI_kR_{k+1})^3)$ flops, which comes at the cost of a squared condition number and possible loss of accuracy. It is possible to construct the matrix of~\eqref{eq:le} without computing any Kronecker product by exploiting the structure of the binary $\mat{S}^{(k)}$ matrices. However, the total runtime of our proposed algorithm is dominated by solving the linear system~\eqref{eq:le} so we will not discuss this particular optimization any further. The key idea of the ALS is to solve~\eqref{eq:le} for varying values of $k$ in a ``sweeping" fashion, starting from the leftmost tensor train core $\ten{A}^{(1)}$ to the rightmost $\ten{A}^{(d)}$ and then back from the rightmost to the leftmost. The numerical stability of the ALS algorithm is guaranteed through an orthogonalization step using the QR decomposition~\cite[p.~A701]{holtz2012alternating}.

\subsection{{Tensor Train Completion Algorithm}}
\label{subsec:ALS}
The pseudocode of the 
proposed tensor train completion (TTC) algorithm is given as Algorithm~\ref{alg:TTC}. First, $d$ binary input matrices $\mat{S}^{(1)}$,$\mat{S}^{(2)}$,\ldots,$\mat{S}^{(d)}$ are constructed as specified in Section~\ref{sec:TNCM}. The tensor train with specified TT-ranks is then initialized according to Section~\ref{subsec:Init}. The ALS is then applied to update the tensor train cores one by one in a sweeping fashion by solving \eqref{eq:le} repeatedly for different values of $k$. Since the tensor train is in site-$d$-mixed-canonical form, the updating starts with $\ten{A}^{(d)}$. One can update the tensor train cores for a fixed amount of sweeps or until the residual falls below a certain threshold. Numerical stability and convergence is guaranteed by the QR factorization step in line $7$. In lines $8$ to 10 of the algorithm the updated tensor $\ten{A}^{(k)}$ is replaced by a reshaping of the orthogonal $\mat{Q}$ matrix and the $\mat{R}$ factor is ``absorbed'' into the next core $\ten{A}^{(k-1)}$. In this way, the resulting tensor train is brought into site-$(k-1)$-mixed-canonical form, before updating $\ten{A}^{(k-1)}$. Global convergence to the solution with unique minimal norm is not guaranteed. Once $\ten{A}^{(2)}$ has been updated, one iteration has completed and the sweep reverses direction with updating the tensor train cores from $\ten{A}^{(1)}$ up to $\ten{A}^{(d-1)}$. Each for-loop in Algorithm~\ref{alg:TTC} therefore corresponds with one iteration.
Depending on the application one can either choose to keep the result in tensor train form or compute the full tensor $\ten{A}$ by summing the tensor train over its auxiliary indices. 
\begin{alg}Tensor completion in tensor train form (TTC)\\
\label{alg:TTC}
\textit{\textbf{Input}}: $d$-way multi-indices of $N$ known entries and their corresponding values $\mat{y}(1),\ldots,\mat{y}(N)$, TT ranks $R_2,\ldots,R_{d}$.\\
\textit{\textbf{Output}}: Completed tensor $\ten{A}$ in tensor train form $\ten{A}^{(1)},\ldots,\ten{A}^{(d)}$.
\begin{algorithmic}[1]
\State Construct $\mat{S}^{(1)}$,$\mat{S}^{(2)}$,\ldots,$\mat{S}^{(d)}$ as specified in Section~\ref{subsec:basicidea}.
\State Initialize the tensor train as specified in Section~\ref{subsec:Init}.
\While {stopping criteria not satisfied } 
\For{k=d,\ldots,2}
\State $\textrm{vec}(\ten{A}^{(k)}) \leftarrow$ solve equation~\eqref{eq:le},~\eqref{eq:leTV} or~\eqref{eq:leTik}
\State $\mat{A}_k \leftarrow$ reshape($\ten{A}^{(k)},[R_{k},I_kR_{k+1}]$)
\State $[\mat{Q},\mat{R}] \gets $ thin QR decomposition of $\mat{A}_k^T$
\State $\ten{A}^{(k)} \leftarrow$ reshape($\mat{Q}^T,[R_{k},I_k,R_{k+1}]$)
\State $\mat{A}_{k-1} \leftarrow$ reshape($\ten{A}^{(k-1)},[R_{k-1}I_{k-1},R_{k}]$)
\State $\ten{A}^{(k-1)} \leftarrow$ reshape($\mat{A}_{k-1}\mat{R}^T,[R_{k-1},I_{k-1},R_{k}]$)
\EndFor
\For{$k=1,\ldots,d-1$}
\State $\textrm{vec}(\ten{A}^{(k)}) \leftarrow$ solve equation~\eqref{eq:le},~\eqref{eq:leTV} or~\eqref{eq:leTik}
\State $\mat{A}_k \leftarrow$ reshape($\ten{A}^{(k)},[R_{k}I_k,R_{k+1}]$)
\State $[\mat{Q},\mat{R}] \gets $ thin QR decomposition of $\mat{A}_k$
\State $\ten{A}^{(k)} \leftarrow$ reshape($\mat{Q},[R_{k},I_k,R_{k+1}]$)
\State $\mat{A}_{k+1} \leftarrow$ reshape($\ten{A}^{(k+1)},[R_{k+1},I_{k+1}R_{k+2}]$)
\State $\ten{A}^{(k+1)} \leftarrow$ reshape($\mat{R}\mat{A}_{k+1},[R_{k+1},I_{k+1},R_{K+2}]$)
\EndFor
\EndWhile
\end{algorithmic}
\end{alg}
The most computationally expensive steps in Algorithm~\ref{alg:TTC} are solving the linear systems in line $5$ and line $13$.

\subsection{Choosing a Dimension Factorization and TT-ranks}
\label{subsec:factordimensions&TTranks}
Color images and videos are 3-way and 4-way tensors, respectively. Converting these tensors directly into tensor trains would then result in 3 or 4 TT-cores with relatively large $I_k$ dimensions. This potentially has a detrimental effect on the runtime of Algorithm~\ref{alg:TTC} as the computational complexity of solving~\eqref{eq:le} is $O((R_kI_kR_{k+1})^3)$. 
The runtime of computing inverse matrices of different sizes varies from one computer to another (the coefficients of the cubic complexity function are dependent on the computer specifications). For the desktop computer used in our experiments, a sharp increase in runtime for matrix inversion was observed when $R_kI_kR_{k+1}\approx 4000$. The problem size was therefore limited to $R_kI_kR_{k+1}\leq 4000$ in our experiments to guarantee a fast completion. This problem size limitation was implemented by factorizing each of the $I_k$ dimensions of the desired completed tensor into more manageable sizes. This implies that each linear index of $I_k$ is split into a corresponding multi-index. This factorization comes at the cost of introducing more TT-cores, resulting in a trade-off between the total number of TT-cores and the total number of parameters per TT-core. The maximal amount of TT-cores is determined by the factorization of each tensor dimension $I_k$ into its prime factors. Specifically, we choose the dimensions of each tensor train core $I_k\leq 10$ combined with TT-ranks $R_kR_{k+1}\leq 400$ for $k=1,\ldots,d$. We will show in an experiment that as long as we follow the above rules, the completion performance of different factorization sizes are similar. Moreover, in case the dimensions are large prime numbers, one can always append zeros to the tensor such that the dimensions can be factored.

\begin{example}
\label{ex:factor}
Suppose we have a $360 \times 640 \times 144 \times 3$ tensor of a color video that consists of 144 frames. The prime factorizations of 360, 640 and 144 are $2^3\times 3^2 \times 5$, $2^7 \times 5$ and $2^4 \times 3^2$, respectively. Separating each of these factors into TT-cores would result in a tensor train of 21 cores. 
The number of cores can, for example, be set to 11 by using the factorization $9 \times 8\times 5\times 4\times 4\times 5\times 8\times 4\times 6\times 6\times 3$. 
\end{example}

Choosing good values for the $d-1$ TT-ranks can be quite tedious when $d$ is large. We therefore propose to choose the values of $R_2$, $R_{\textrm{mid}}$ (``mid" stands for middle) and $R_d$ and automatically determine the remaining TT-ranks $R_{k+1} (2\leq k \leq d-2)$ as $\textrm{min}(R_kI_k,R_{\textrm{mid}})$. In this way, a uniform plateau of TT-ranks equal to $R_{\textrm{mid}}$ is obtained. The TT-ranks $R_2,R_d$ need to be chosen while keeping in mind that \mbox{$R_2 \leq I_1$} and $R_d \leq I_d$. When the tensor represents either a color image or color video, the last dimension will typically be 3 for the color channels. In this case, we always set $R_d=3$ and additionally choose the value of $R_{d-1}$.
\begin{example}
We revisit the $360 \times 640 \times 144 \times 3$ tensor of Example~\ref{ex:factor}, together with the 11-core dimension factorization $9 \times 8\times 5\times 4\times 4\times 5\times 8\times 4\times 6\times 6\times 3$. Choosing $R_2=5, R_{\textrm{mid}}=5$, $R_{d-1}=R_{10}=5$ and $R_{11}=3$ then results in $R_3=\cdots =R_{9}=5$.
\end{example}


\subsection{Total Variation Regularization}
\label{subsec:TV}
The low TT-rank constraint in~\eqref{eqn:ude} can be interpreted as a global feature, as it pertains to the construction of the whole tensor. Better completion results can be obtained from the addition of local constraints on the completed tensor entries. The notion of local smoothness, as described by TV, is such a local feature, and is particularly useful when the tensor represents visual data. The addition of TV terms to~\eqref{eqn:ude} is quite straightforward, resulting in
\begin{align}
\label{eqn:TVprob}
\nonumber \min_{\ten{A}\in\mathcal{S}^{(d)}_{\text{TT}}} \;||\mat{S}^T \;\textrm{vec}(\ten{A}) - \mat{y}||_2^2 + &\sum_{p=1}^{P} \lambda_p\, ||\ten{A} \times_p \mat{D}_p ||_F^2,\\
\textrm{\textrm{s.t.}} \;\ten{A}\in\mathcal{S}^{(d)}_{\text{TT}}\text{, and TT-rank}(\ten{A}&) = (R_1,R_2,\ldots,R_{d}),
\end{align}
where $\mat{D}_p \in \mathbb{R}^{(I_p-1) \times I_p}$ has entries \mbox{$\mat{D}_p(i,i)=1$} and \mbox{$\mat{D}_p(i,i+1)=-1$}. For notational convenience we make $\mat{D}_p$ square by appending a bottom row of zeros. Minimizing each of the $||\ten{A} \times_p \mat{D}_p ||_F^2$ terms ensures that the completed tensor entries do not differ too much from their neighbors along modes $1,\ldots,P$, which encodes local smoothness. For image and video tensors we have $P=2$, as the smoothness is only required in the width and length dimensions. Solving~\eqref{eqn:TVprob} requires the modification of~\eqref{eq:le} with additional matrix terms. To derive how exactly~\eqref{eq:le} needs to be modified, we will first ignore the fact that we can factorize the dimensions of the original tensor as discussed in Section~\ref{subsec:factordimensions&TTranks}. It is worthwhile to stress that to the best of our knowledge, this is the first time that the Total Variation regularization term is exploited and fully incorporated in tensor train form.

\begin{figure}[tb]
\begin{center}
\input{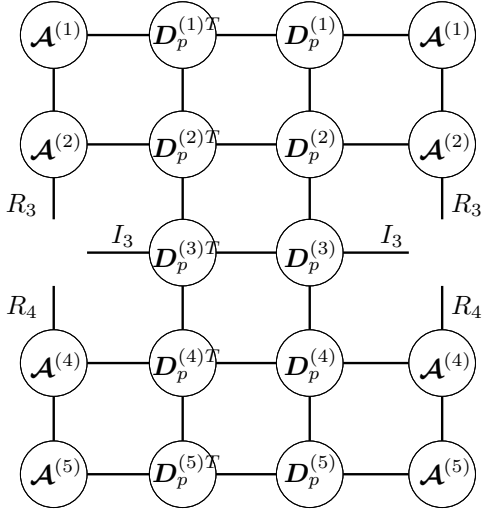}
\caption{Tensor network diagram of $\mat{W}_p^T\mat{W}_p$ when updating $\ten{A}^{(3)}$ in a tensor train of 5 TT-cores.}
\label{fig:WtW}
\end{center}
\end{figure}

The derivation of this modification is very similar to the traditional ALS derivation, where the input matrices $\mat{S}^{(k)}$ are now replaced with the matrices
\begin{align*}
\mat{D}^{(k)}_p :=\left\{ 
\begin{aligned}
& \mat{D}_p,    &  & \text{if} \  k=p, \\ 
& \mat{I}_{I_k},  &  & \text{otherwise} \ . 
\end{aligned}
\right.
\end{align*}
Indeed, the $p$-th TV regularization term can now be rewritten as $(\mat{D}^{(d)}_p \otimes \cdots \otimes \mat{D}^{(2)}_p \otimes \mat{D}^{(1)}_p)\, \textrm{vec}(\ten{A})$. The only difference between the TV term and $\mat{S}^T\,\textrm{vec}(\ten{A})$ is that $\mat{S}$ consists of a Khatri-Rao product, while the TV term contains a Kronecker product. Just as in Section~\ref{subsec:ALS} we consider the contractions of the TT-cores of $\ten{A}$ with these new ``input'' matrices
\begin{align*}
\mat{W}_{<k,p} &:= (\ten{A}^{(1)} \times_2 \mat{D}^{(1)}_p)  \ldots (\ten{A}^{(k-1)} \times_2 \mat{D}^{(k-1)}_p),\\
\mat{W}_{>k,p} &:= (\ten{A}^{(k+1)} \times_2 \mat{D}^{(k+1)}_p ) \ldots (\ten{A}^{(d)} \times_2 \mat{D}^{(d)}_p),
\end{align*}
which allows us to define the matrix
\begin{align*}
\mat{W}_p &:= (\mat{W}_{>k,p} \otimes \mat{D}^{(k)}_p \otimes \mat{W}_{<k,p}) \in \mathbb{R}^{\prod_{i=1}^d I_i \times R_kI_kR_{k+1}}.
\end{align*}
Denoting the matrix in~\eqref{eq:le} by $\mat{B}$, the modified linear system is then
\begin{align}
\label{eq:leTV}
(\mat{B}^T\mat{B} + \sum_{p=1}^P \lambda_p \mat{W}_p^T\mat{W}_p)\; \textrm{vec}(\ten{A}^{(k)}) &= \mat{B}^T\,\mat{y}.
\end{align}
The number of rows of $\mat{W}_p$ grows exponentially with $d$. Fortunately, it is not necessary to explicitly construct this matrix as it is possible to compute $\mat{W}_p^T\mat{W}_p$ directly from the tensor network diagram depicted in Figure~\ref{fig:WtW} for the case $d=5, k=3$. This computation can be done efficiently by exploiting the fact that most of the $\mat{D}_p^{(i)}$ matrices are the unit matrix. In addition, the site-$k$-mixed-canonical form of the tensor train of $\ten{A}$ can also be exploited. Suppose for example that $p=1$ in Figure~\ref{fig:WtW}. The site-3-mixed-canonical form of the tensor train of $\ten{A}$ then implies that  both $\ten{A}^{(4)}$ and $\ten{A}^{(5)}$ are right-orthogonal. Since $p=1$, both $\mat{D}_1^{(4)}$ and $\mat{D}_1^{(5)}$ matrices are unit matrices. The summation of the two bottom rows in Figure~\ref{fig:WtW} then results in an $R_4 \times R_4$ unit matrix and can be skipped. Also note that all vertical edges between the different $\mat{D}_p^{(k)}$ matrices have dimensions 1.

When the dimensions of the original tensor are factorized as discussed in Section~\ref{subsec:factordimensions&TTranks}, then only one minor modification is required. The TT-SVD algorithm needs to be applied to $\mat{D}_p$ to transform this matrix into a tensor train matrix (TT matrix)~\cite{ttmatrix} according to the same dimension factorization. The corresponding node in the tensor network diagram is then replaced by the corresponding TT matrix. As a result, not all vertical edges will have dimension 1 anymore. It turns out that the TT matrix-ranks of $\mat{D}_p$ are always uniformly 3, irrespective of the factorization of the dimensions or number of cores. This low TT-rank feature of the $\mat{D}_p$ matrix is favorable for the computation of $\mat{W}_p^T\mat{W}_p$ from the tensor network.

\begin{example}
Suppose we have a $1024 \times 1024 \times 3$ color image and we factor each of the 1024 dimensions into $4^5$. The tensor train of $\textrm{vec}(\ten{A})$ hence consists of 11 cores. The $1024 \times 1024$ $\mat{D}_p$ matrix is then converted into a TT matrix of 5 cores with ranks $R_2=\cdots=R_5=3$. When computing the $\mat{W}^T_1\mat{W}_1$ term, we then have that the first five TT-cores of $\ten{A}$ are connected with the TT matrix of $\mat{D}_1$, while the remaining cores are connected to identity matrices. 
\end{example}

\begin{table*}[t]
\centering
\begin{center}
\caption{Performance (RSEs) of nine algorithms on eight benchmark images}
\resizebox{\textwidth}{!}
{
\begin{tabular}{
@{}lrrrrrrrrrrrrrrrr@{}}
\centering
\multirow{2}{*}{Method} & \multicolumn{2}{c}{House} & \multicolumn{2}{c}{River} & \multicolumn{2}{c}{Bridge} & \multicolumn{2}{c}{Man} & \multicolumn{2}{c}{Lena} & \multicolumn{2}{c}{Peppers} & \multicolumn{2}{c}{Baboon} & \multicolumn{2}{c}{Airplane} \\ 
& RSE & Time(s) & RSE & Time(s) & RSE & Time(s) & RSE & Time(s) & RSE & Time(s) & RSE & Time(s) & RSE & Time(s) & RSE & Time(s)\\\hline
TNN & $0.221$ & $44.9$ & $0.193$ & $43.0$ & $0.159$ & $60.7$ & $0.231$ & $33.5$ & $0.221$ & $13.6$ & $0.295$ & $16.6$ & $0.246$ & $15.7$ & $0.168$ & $18.1$\\
HaLRTC & $0.209$ & ${7.9}$ & $0.182$ & $6.3$ & $0.153$ & $15.6$ & $0.221$ & $20.3$ & $0.206$ & $10.6$ & $0.265$ & ${5.2}$ & $0.222$ & $4.1$ & $0.162$ & $8.4$\\
FaLRTC & $0.210$ & $30.0$ & $0.181$ & $10.0$ & $0.153$ & $26.6$ & $0.218$ & $10.6$ & $0.202$ & $5.1$ & $0.253$ & $5.3$ & $0.228$ & $7.0$ & $0.161$ & $12.3$\\
LRTC-TV & ${0.138}$ & $163.1$ & ${0.123}$ & $178.4$ & ${0.126}$ & $165.8$ & ${0.130}$ & $185.7$ & $0.107$ & $97.4$ & ${0.127}$ & $99.3$ & ${0.159}$ & $99.7$ & $0.102$ & $98.9$\\
TMac & $0.226$ & $12.6$ & $0.199$ & $8.4$ & $0.168$ & $11.4$ & $0.237$ & $12.0$ & $0.240$ & $5.6$ & $0.287$ & $7.7$ & $0.236$ & $6.0$ & $0.173$ & $8.2$\\
TMac-TT & $0.225$ & $21.0$ & $0.191$ & $28.6$ & $0.181$ & $28.1$ & $0.203$ & $20.7$ & $0.165$ & $10.2$ & $0.215$ & $9.6$ & $0.201$ & $7.4$ & $0.145$ & $9.8$\\
TR-ALS & $0.170$ & $8.2$ & $0.154$ & $7.6$ & $0.153$ & ${9.5}$ & $0.185$ & $8.5$ & $0.160$ & $5.4$ & $0.198$ & $7.2$ & $0.205$ & $5.1$ & $0.141$ & ${3.6}$\\
TTC & $0.161$ & $20.7$ & $0.153$ & $38.1$ & $0.161$ & $39.9$ & $0.171$ & $27.9$ & $0.158$ & $15.3$ & $0.194$ & $20.5$ & $0.212$ & $29.5$ & $0.142$ & $27.6$\\
TTC-TV & $0.151$ & $23.1$ & $0.136$ & $40.4$ & $0.141$ & $40.9$ & $0.155$ & $29.8$ & $0.131$ & $15.8$ & $0.175$ & $22.2$ & $0.176$ & $31.0$ & $0.111$ & $28.7$\\
\end{tabular}
}
\label{tbl:bm}
\end{center}
\end{table*}  
\begin{figure}[t]
\centering
\includegraphics[width=0.48\textwidth]{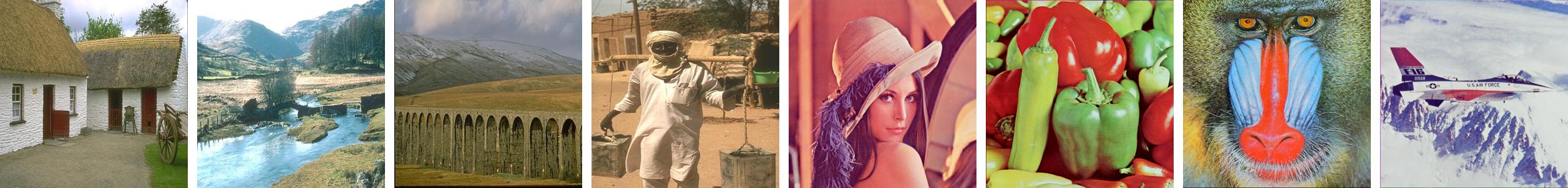}
\caption{Ground-truth of eight small benchmark images.}
\label{fig:groundtruth}
\end{figure}
\begin{figure*}
\centering
\subfigure[Observed]{
\begin{minipage}[b]{0.09\linewidth}
\includegraphics[width=1\linewidth]{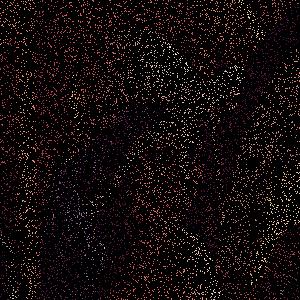}\vspace{4pt}
\includegraphics[width=1\linewidth]{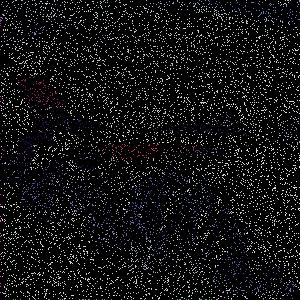}\vspace{4pt}
\end{minipage}}
\subfigure[TNN]{
\begin{minipage}[b]{0.09\linewidth}
\includegraphics[width=1\linewidth]{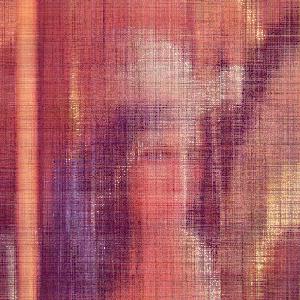}\vspace{4pt}
\includegraphics[width=1\linewidth]{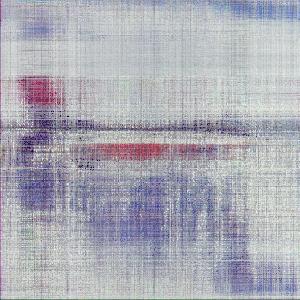}\vspace{4pt}
\end{minipage}}
\subfigure[HaLRTC]{
\begin{minipage}[b]{0.09\linewidth}
\includegraphics[width=1\linewidth]{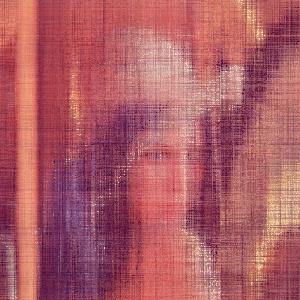}\vspace{4pt}
\includegraphics[width=1\linewidth]{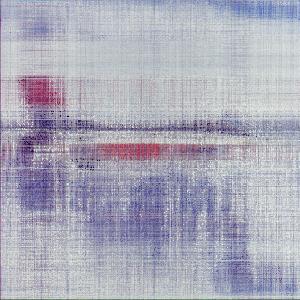}\vspace{4pt}
\end{minipage}}
\subfigure[LRTC-TV]{
\begin{minipage}[b]{0.09\linewidth}
\includegraphics[width=1\linewidth]{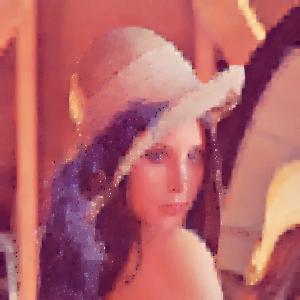}\vspace{4pt}
\includegraphics[width=1\linewidth]{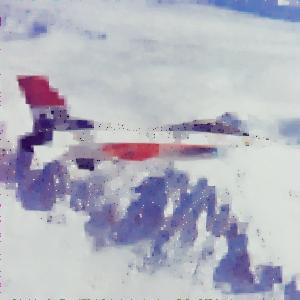}\vspace{4pt}
\end{minipage}}
\subfigure[TMac]{
\begin{minipage}[b]{0.09\linewidth}
\includegraphics[width=1\linewidth]{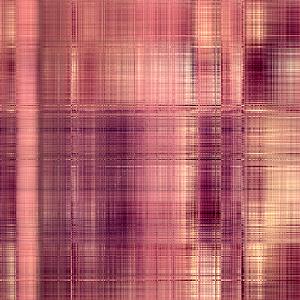}\vspace{4pt}
\includegraphics[width=1\linewidth]{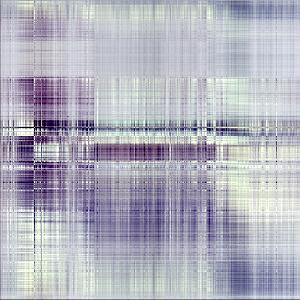}\vspace{4pt}
\end{minipage}}
\subfigure[TMac-TT]{
\begin{minipage}[b]{0.09\linewidth}
\includegraphics[width=1\linewidth]{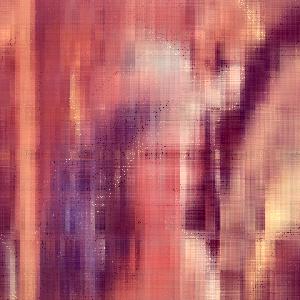}\vspace{4pt}
\includegraphics[width=1\linewidth]{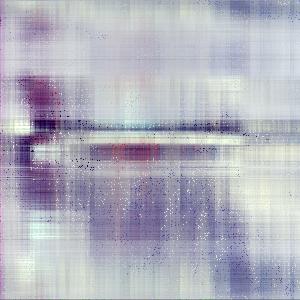}\vspace{4pt}
\end{minipage}}
\subfigure[TR-ALS]{
\begin{minipage}[b]{0.09\linewidth}
\includegraphics[width=1\linewidth]{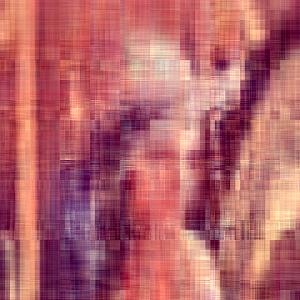}\vspace{4pt}
\includegraphics[width=1\linewidth]{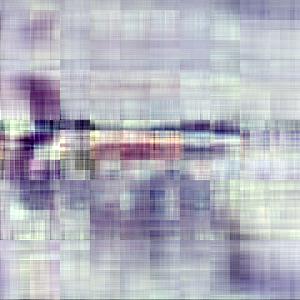}\vspace{4pt}
\end{minipage}}
\subfigure[TTC]{
\begin{minipage}[b]{0.09\linewidth}
\includegraphics[width=1\linewidth]{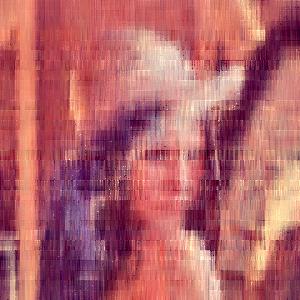}\vspace{4pt}
\includegraphics[width=1\linewidth]{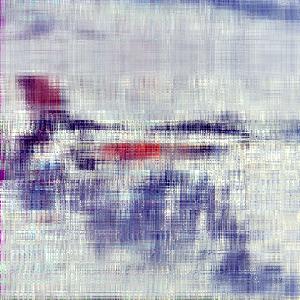}\vspace{4pt}
\end{minipage}}
\subfigure[TTC-TV]{
\begin{minipage}[b]{0.09\linewidth}
\includegraphics[width=1\linewidth]{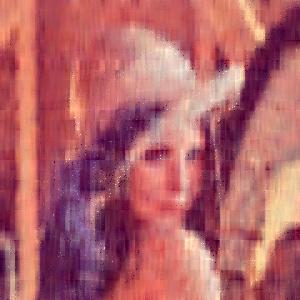}\vspace{4pt}
\includegraphics[width=1\linewidth]{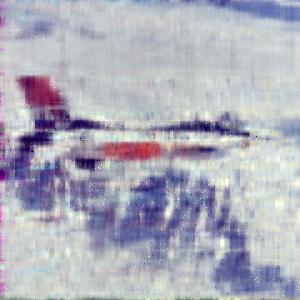}\vspace{4pt}
\end{minipage}}
\caption{Image inpaintings of \textit{Lena}, and \textit{Airplane} benchmarks by seven algorithms.}
\label{fig:visual}
\end{figure*}

The values of the $\lambda_p$ parameters were fixed to $0.5$ in~\cite{li2017low}. Through experiments, we found that choosing values between $0$ and $10$ for the $\lambda_p$ parameters resulted in consistently improved completion results compared to the standard TTC ALS method. Specifically, an initial $\lambda_p$ of $1$ is chosen throughout our later experiments. It is also possible to adjust the values of the $\lambda_p$ parameters during the execution of Algorithm~\ref{alg:TTC}. As the estimate of the completed image improves over the iterations, it might not be necessary to keep enforcing the local smoothness. The heuristic we propose is to multiply the $\lambda_p$ parameters with the relative error on the observed errors at the end of each iteration and use those values for the next iteration. In practice, this results in a more monotonic convergence of the relative error as a function of the iterations.

\subsection{Tikhonov regularization}
\label{subsec:Tikhonov}
In addition to the TV terms, Tikhonov regularization has also been considered in~\cite[p.~2213]{li2017low}. The tensor completion problem is then written as the following optimization problem
\begin{align}
\label{eqn:Tikprob}
\nonumber \min_{\ten{A}\in\mathcal{S}^{(d)}_{\text{TT}}} \,||\mat{S}^T \textrm{vec}(\ten{A}) - \mat{y}||_2^2 + \sum_{p=1}^{P} \lambda_p& ||\ten{A} \times_p \mat{D}_p ||_F^2 + \gamma || \ten{A}||_F^2,\\
\textrm{\textrm{s.t.}} \;\ten{A}\in\mathcal{S}^{(d)}_{\text{TT}}\text{, and TT-rank}(\ten{A}&) = (R_1,R_2,\ldots,R_{d}).
\end{align}
Using the fact that the tensor train is in site-$k$-mixed-canonical form together with~\eqref{eqn:norm}, the update step in the ALS algorithm is then modified to
\begin{align}
\label{eq:leTik}
(\mat{B}^T\mat{B} + \sum_{p=1}^P \lambda_p \mat{W}_p^T\mat{W}_p + \gamma \,\mat{I})\; \textrm{vec}(\ten{A}^{(k)}) &= \mat{B}^T\mat{y},
\end{align}
where $\mat{I}$ is an identity matrix of size $R_kI_kR_{k+1}$. The addition of a Tikhonov regularization term therefore comes at zero additional computational cost.

\section{{Experimental Results}}
\label{sec:E}

In this section our proposed algorithm is compared extensively with state-of-the-art tensor completion methods in three experiments. The scalability of our algorithm in particular is demonstrated through the last two experiments. All experiments were run in MATLAB2018a on a desktop computer with an Intel i5 quad-core processor at 3.2GHz and 16GB RAM. A MATLAB implementation of Algorithm~\ref{alg:TTC}, together with necessary data files for reproducing all experimental results can be downloaded from~\url{https://github.com/IRENEKO/TTC}. The values of all tuning parameters used in these experiments as well as benchmark images are all given in the supplemental materials. First, we apply Algorithm~\ref{alg:TTC} to complete eight color images with approximate size $300\times 300$ and compare its efficacy in runtime and completion accuracy with seven other state-of-the-art tensor completion algorithms. The two best methods from the first experiment are then compared with our algorithm in the second experiment for the completion of three color images with approximate size $4000 \times 4000$. The increased dimensions allow us to highlight the scalability of Algorithm~\ref{alg:TTC} with these state-of-the-art methods. Another way to assess the scalability of our method is to apply it on higher order tensors. We therefore also compare Algorithm~\ref{alg:TTC} with the two best methods from Experiment 1 to complete a color video.
The completion accuracy of all methods is measured either by the relative standard error (RSE)
\begin{align*}
\textrm{RSE}&=\frac{||\ten{A}-\ten{\hat{A}}||_F}{||\ten{A}||_F},
\end{align*} 
or the peak signal-to-noise ratio (PSNR)
\begin{align*}
\textrm{PSNR}&=20\;\log_{10}(\textrm{MAX}_\textrm{I})-10\;\log_{10}(\textrm{MSE}),
\end{align*}
where $\ten{\hat{A}}$ is the completed tensor, $\textrm{MAX}_\textrm{I}$ is the maximum possible pixel value and \textrm{MSE} is the mean square error $||\ten{A}-\ten{\hat{A}}||_F/\textrm{numel}(\ten{A})$, where $\textrm{numel}(\ten{A})$ denotes the total number of entries in $\ten{A}$.

\subsection{{Small Image Inpainting}}

\begin{table}[t]
\scriptsize
\begin{center}
\caption{Experiment 1 dimension settings.}
\begin{tabular}{@{}lcc@{}}
Image & Resized dimensions & Dimension factorization\\
\hline
House, River\\
Bridge, Man & $324 \times 486 \times 3$ & $9 \times 6 \times 6 \times 6 \times 9 \times 9 \times 3$\\
Lena, Peppers\\
Baboon, Airplane & $300 \times 300 \times 3$ & $10 \times 6 \times 5 \times 5 \times 6 \times 10 \times 3$
\end{tabular}
\label{tbl:pre1}
\end{center}
\end{table}  

Eight benchmark images, shown in Figure~\ref{fig:groundtruth}, from the Berkeley Segmentation dataset\footnote{https://www.eecs.berkeley.edu/Research/Projects/CS/vision/bsds/} and USC-SIPI image database\footnote{http://sipi.usc.edu/database/database.php} were used to compare the performance of our proposed method with state-of-the-art completion algorithms in terms of completion accuracy and runtime. TTC denotes Algorithm~\ref{alg:TTC} without any TV or Tikhonov regularization and TTC-TV denotes Algorithm~\ref{alg:TTC} with TV regularization. Tikhonov regularization did not improve the results significantly and is therefore not considered.
Table~\ref{tbl:pre1} lists the dimensions of the benchmark images and the dimension factorizations used for the tensor train methods. All images were resized using bicubic interpolation through the MATLAB ``imresize'' command. The eight images are grouped into two groups with slightly different dimensions, which is a first attempt to determine how sensitive the runtime and completion accuracy of all methods are with respect to dimension size.
Only $10\%$ of each benchmark image was retained, whereby for each missing pixel all color information was removed. 
We \textbf{fine-tune} the hyper-parameters of the seven algorithms in comparison to give the best RSE scores on the images. 

For the proposed TTC method, the TT-ranks are determined by \textbf{cross-validation}~\cite{li2016structured} on the completion error of a held-out $1\%$ entries. That is, for an image with only $10\%$ observed entries, $1/10$ of these known entries are kept for validation. We perform $10$ trials in each cross-validation experiment and the TT-ranks that give the lowest average error of held-out entries are chosen. The same TT-ranks are shared by TTC-TV experiments for the same image.

\subsubsection{State-of-the-art methods}
Algorithm~\ref{alg:TTC} is compared with the following state-of-the-art methods: TNN\footnote{http://www.ece.tufts.edu/~shuchin/software.html}, HaLRTC, FaLRTC\footnote{http://www.cs.rochester.edu/u/jliu/code/TensorCompletion.zip}, LRTC-TV\footnote{https://xutaoli.weebly.com/}, 
TMac\footnote{http://www.math.ucla.edu/~wotaoyin/papers/codes/TMac.zip}, TMac-TT\footnote{https://sites.google.com/site/jbengua/home/projects/efficient-tensor-completion-for-color-image-and-video-recovery-low-rank-tensor-train}, and TR-ALS\footnote{https://github.com/wangwenqi1990/TensorRingCompletion}. These methods represent four different approaches towards the completion problem. The TNN method aims at minimizing the number of nonzeros in the tensor multi-rank, which is later relaxed to minimize the nuclear norm of a matrix constructed by frontal slices of the three way tensor~\cite{zhang2014novel}. HaLRTC, FaLRTC~\cite{liu2013tensor} and LRTC-TV~\cite{li2017low} are on the other hand minimizing the sum of nuclear norms of the unfolded matrices. The TMac algorithm~\cite{xu2015parallel} includes two different schemes, TMac-inc and TMac-dec, depending on different rank adjustment strategies. Here we only compare with TMac-inc as it shows a better performance than TMac-dec in our experiments. The extension of TMac algorithm TMac-TT~\cite{bengua2017efficient} is also considered. The TR-ALS algorithm~\cite{wang2017efficient} also uses a tensor train of $\ten{A}$ but with $R_1>1$. Moreover, it employs a different ALS for updating each tensor train core wherein each slice of the core is updated sequentially, unlike our way of updating the whole core at once. In these benchmarking algorithms, all MATLAB implementations written by the original authors were used. 

\begin{table}
\scriptsize
\begin{center}
\caption{Average runtime over the eight images of each algorithm to obtain its best RSE score and the corresponding average runtime of the TTC algorithm to obtain the same RSE score.}
\begin{tabular}{@{}lrrrrrrr@{}}
& TNN & HaLRTC & FaLRTC 
& TMac & TMac-TT & TR-ALS\\\hline
Average \\
runtime (s) & $30.8$ & $9.8$ & $13.4$ 
& $9.0$ & $17.2$ & $6.9$\\
Average TTC\\
runtime (s) & $1.0$ & $1.6$ & $1.6$ 
& $0.8$ & $3.0$ & $4.3$ \\ [1.ex]
Speedup & $30.8\times $ & $6.1\times $ & $8.4\times $ 
& $11.3\times $ & $5.7\times $ & $1.6\times $\\
\end{tabular}
\label{tbl:efficient}
\end{center}
\end{table}  

\begin{table*}[t]
\centering
\begin{center}
\caption{Performance (RSEs) of nine algorithms on eight benchmark images with $1\%$ observed entries}
\resizebox{\textwidth}{!}
{
\begin{tabular}{
@{}lrrrrrrrrrrrrrrrr@{}}
\centering
\multirow{2}{*}{Method} & \multicolumn{2}{c}{House} & \multicolumn{2}{c}{River} & \multicolumn{2}{c}{Bridge} & \multicolumn{2}{c}{Man} & \multicolumn{2}{c}{Lena} & \multicolumn{2}{c}{Peppers} & \multicolumn{2}{c}{Baboon} & \multicolumn{2}{c}{Airplane} \\ 
& RSE & Time(s) & RSE & Time(s) & RSE & Time(s) & RSE & Time(s) & RSE & Time(s) & RSE & Time(s) & RSE & Time(s) & RSE & Time(s)\\\hline
TNN & $0.523$ & $11.0$ & $0.469$ & $12.4$ & $0.431$ & $55.3$ & $0.498$ & $31.4$ & $0.842$ & $2.3$ & $0.662$ & $6.0$ & $0.585$ & $12.1$ & $0.500$ & $15.6$\\
LRTC-TV & $0.228$ & $154.5$ & $0.191$ & $150.8$ & $0.202$ & $138.3$ & $0.296$ & $150.1$ & $0.324$ & $81.1$ & $0.453$ & $77.6$ & $0.353$ & $77.0$ & $0.198$ & $71.5$\\
TMac-TT & $0.360$ & $19.4$ & $0.278$ & $18.7$ & $0.241$ & $17.8$ & $0.326$ & $17.1$ & $0.301$ & $9.6$ & $0.409$ & $20.9$ & $0.292$ & $11.5$ & $0.209$ & $9.3$\\
TR-ALS & $0.350$ & $2.2$ & $0.260$ & $3.0$ & $0.251$ & $2.8$ & $0.322$ & $2.3$ & $0.347$ & $1.8$ & $0.457$ & $1.6$ & $0.339$ & $1.3$ & $0.244$ & $1.3$\\
TTC-TV & $0.210$ & $3.2$ & $0.197$ & $5.7$ & $0.199$ & $1.4$ & $0.232$ & $2.4$ & $0.232$ & $0.9$ & $0.295$ & $0.7$ & $0.250$ & $1.1$ & $0.176$ & $1.0$\\
\end{tabular}
}
\label{tbl:bm2}
\end{center}
\end{table*}  
\begin{figure*}[t]
\centering
\subfigure[LRTC-TV: $\textrm{PSNR}=24.577\textrm{dB}$, $\textrm{RSE}=0.107$.]{
\begin{minipage}[b]{0.4\linewidth}
\includegraphics[width=1\linewidth]{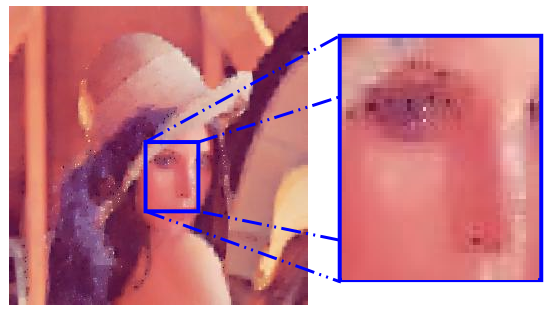}\vspace{4pt}
\end{minipage}}
\subfigure[TTC-TV: $\textrm{PSNR}=24.423\textrm{dB}$, $\textrm{RSE}=0.109$.]{
\begin{minipage}[b]{0.4\linewidth}
\includegraphics[width=1\linewidth]{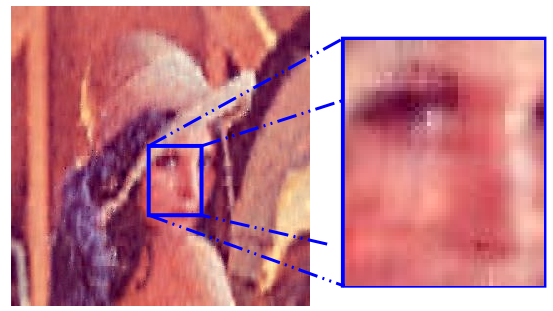}\vspace{4pt}
\end{minipage}}
\caption{Image inpaintings of \textit{Lena} by LRTC-TV and TTC-TV. LRTC-TV scores better on the PSNR and RSE metrics but lacks details due to oversmoothing.}
\label{fig:Lena}
\end{figure*}

\subsubsection{Comparison with state-of-the-art}

The completion accuracy measured as the RSE of all algorithms and their corresponding runtimes are reported in Table~\ref{tbl:bm}, while the completed images for \textit{Lena} and \textit{Airplane} are shown in Figure~\ref{fig:visual}. The completed images obtained by the FaLRTC is indistinguishable from HaLRTC and has therefore been omitted from Figure~\ref{fig:visual}. The proposed TTC-TV algorithm with TTC cross-validated TT-ranks outperforms the other six algorithms, excluding LRTC-TV, in terms of RSE for all eight benchmark images at the cost of an overall larger runtime than TR-ALS. Moreover, although TTC only reaches lower RSEs compared to its competitors on five out of eight benchmark images, we remark that it is able to obtain the same or even lower RSEs with less runtime with fine-tuned TT-ranks.
Table~\ref{tbl:efficient} lists the runtime for each method averaged over all images to obtain its best RSE, together with the average runtime over all images for the TTC algorithm to obtain an identical RSE. The TTC algorithm is seen to be faster than the state-of-the-art methods, with average speedups of at least 10 compared with TMac and TNN. Moreover, it is worth noting that the proposed TTC-TV algorithm results in consistently smaller RSEs compared with those achieved by TTC at similar runtimes. These better RSE values obtained with TTC-TV result in completed images that are smoother while still preserving details. LRTC-TV, an extension of the LRTC methods with total variation terms, reaches the lowest RSEs in all benchmark images at the cost of taking $4.7$ times longer than TTC-TV on average. 

\subsubsection{RSE is not enough}
Table~\ref{tbl:bm} seems to indicate that LRTC-TV consistently obtains better completion results over all other methods. However, we argue that it is sometimes not enough to evaluate the completed images by their corresponding RSEs. A visual inspection of the completed images still remains the best way to compare results. We illustrate this point by comparing the completed \textit{Lena} image using both the LRTC-TV and TTC-TV methods in Figure~\ref{fig:Lena}. Both the TT-rank and the number of iterations used in TTC-TV were increased to obtain an RSE score that is quite close to the LRTC-TV method. Although the LRTC-TV method has better RSE and PSNR scores, one can see that detailed features such as the eyes and feathers are blurred by the LRTC-TV method due to oversmoothing. This oversmoothing was observed for any nonzero value of the $\lambda$ tuning parameters.

\subsubsection{Influence of smaller portions of observed entries}
In what follows, we explore the influence of smaller percentages of observed entries on the completion accuracy and runtime by recovering images with $1\%$ observed. As before, we fine-tune the tuning parameters of different algorithms except the proposed TTC-TV, for which cross-validated hyper-parameters are used. LRTC-TV, TMac-TT and TTC-TV generally perform better than LRTCs, TMac and TTC, respectively, according to Table~\ref{tbl:bm}, hence we omit the latter three in Table~\ref{tbl:bm2}. Besides the proposed TTC-TV method, LRTC-TV outperforms the TNN, TMac-TT and TR-ALS algorithms in five out of eight benchmark images at the cost of $6, 7, 55$ times longer average runtime, respectively. TMac-TT reaches lower RSEs than LRTC-TV when completing the remaining three images. Moreover, the proposed TTC-TV outperforms all other methods in all benchmark images except \textit{River}, where a similar RSE as that obtained by LRTC-TV is reached. It is also worthwhile to note that TTC-TV demonstrates a speedup of up to $55$ times compared with LRTC-TV in the above investigation.

\subsubsection{Effect of TT-Initialization on ALS convergence}
\begin{figure}[t]
\centering
\includegraphics[width=0.5\textwidth]{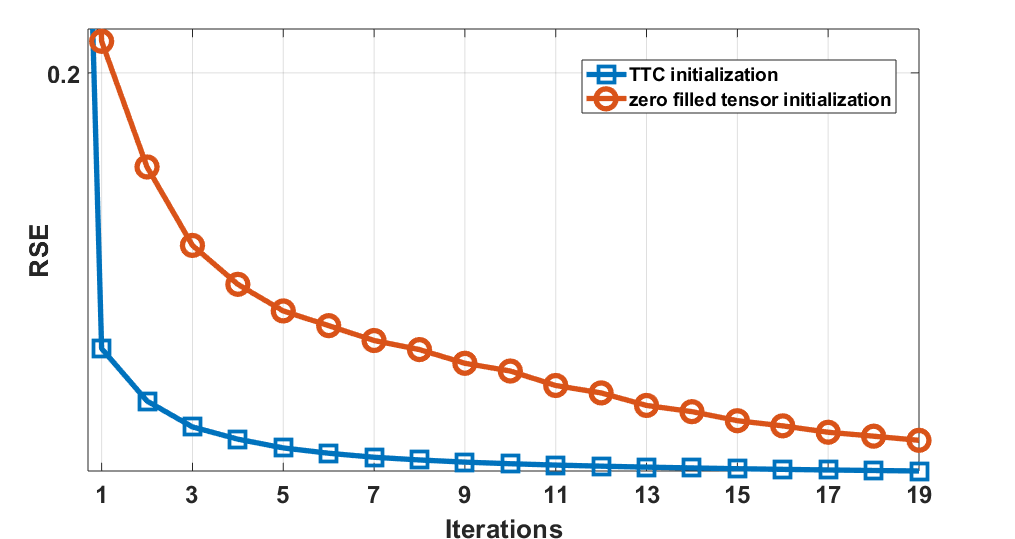}
\caption{The RSE as a function of the iteration count for two different initialization methods of the TTC algorithm.}
\label{fig:init}
\end{figure}

We also investigated whether our proposed initialization method results in better convergence behavior of ALS. An alternative initialization method, called the ``missing entry zero filled'' tensor initialization, is described in~\cite{wang2017efficient} and used in the TR-ALS method. Figure~\ref{fig:init} depicts the RSEs of the completed \textit{Man} benchmark image as a function of the iteration count using these two different initializations. If we assume that the RSE obtained after 8 iterations when using the proposed initialization method is taken as a threshold to stop the TTC algorithm, then the ``missing entry zero-filled" tensor initialization would still not have converged after 20 iterations. In general, using our proposed initialization method always resulted in the RSE curve tapering off very fast over iterations. However, modifying the code for TNN, LRTC-TV, TMac, TMac-TT to the proposed initialization method is not trivial. For HaLRTC, FaLRTC and TR-ALS, the algorithms get stuck in local optima and achieve high RSEs ($\approx 0.9$) with the interpolation initialization.

\subsubsection{Influence of different factorization sizes}
\label{subsubsec:factorization}
In Section~\ref{subsec:factordimensions&TTranks}, a general guideline for choosing a specific factorization size of a tensor is provided. Here we use the \textit{House} benchmark image as an example to show the influence of different factorization sizes on both the obtained RSE, total runtime and problem size ($R_kI_kR_{k+1}$), all listed in Table~\ref{tbl:factorization}. The total number of iterations of Algorithm~\ref{alg:TTC} is fixed for all factorizations. As shown in the table, similar RSEs are obtained when the maximal dimension in the factorization is limited to 9. A slight improvement of the runtime is observed when the number of tensor train cores is reduced from 13 down to 7 due to the trade-off between the number of cores and problem size as discussed in Section~\ref{subsec:factordimensions&TTranks}. Both  the RSE and total runtime are seen to degrade as the dimensions in the factorization are further increased, due to the fast growing problem size. These observations are consistent with the rules specified in Section~\ref{subsec:factordimensions&TTranks}.
\begin{table}[t]
\scriptsize
\begin{center}
\caption{Influence of different factorization sizes.}
\begin{tabular}{@{}lrrr@{}}
Dimension factorization & RSE & Time(s) & Problem size\\
\hline
$324 \times 486 \times 3$ & $0.200$ & $820.4$ & $14580$\\
$18^2 \times 18 \times 27 \times 3$ & $0.174$ & $28.4$ & $2160$\\
$6^2 \times 9 \times 6 \times 9^2 \times 3$ & $0.161$ & $10.4$ & $1683$\\
$2^2 \times 3^4 \times 2 \times 3^5 \times 3$ & $0.161$ & $14.4$ & $1083$
\end{tabular}
\label{tbl:factorization}
\end{center}
\end{table}  

\subsection{Large Image Inpainting}

\begin{table}[t]
\scriptsize
\begin{center}
\caption{Experiment 2 dimension settings.}
\begin{tabular}{@{}lrrr@{}}
Image & Original dimensions & Dimension factorization \\
\hline
Dolphin & $4000 \times 3000 \times 3$ & $2^5 \times 5^3 \times 2^3 \times 3 \times 5^3 \times 3$ \\
Water Nature Fall & $6000 \times 4000 \times 3$ & $2^4 \times 3 \times 5^3 \times 2^5 \times 5^3  \times 3$\\
Orion nebula  & $6000 \times 6000 \times 3$ & $2^4 \times 3 \times 5^3 \times 2^4 \times 3 \times 5^3 \times 3$
\end{tabular}
\label{tbl:pre2}
\end{center}
\end{table}  

\begin{figure*}[tbh]
\centering
\subfigure[Dolphin.]{
\begin{minipage}[b]{0.45\linewidth}
\includegraphics[width=1\linewidth]{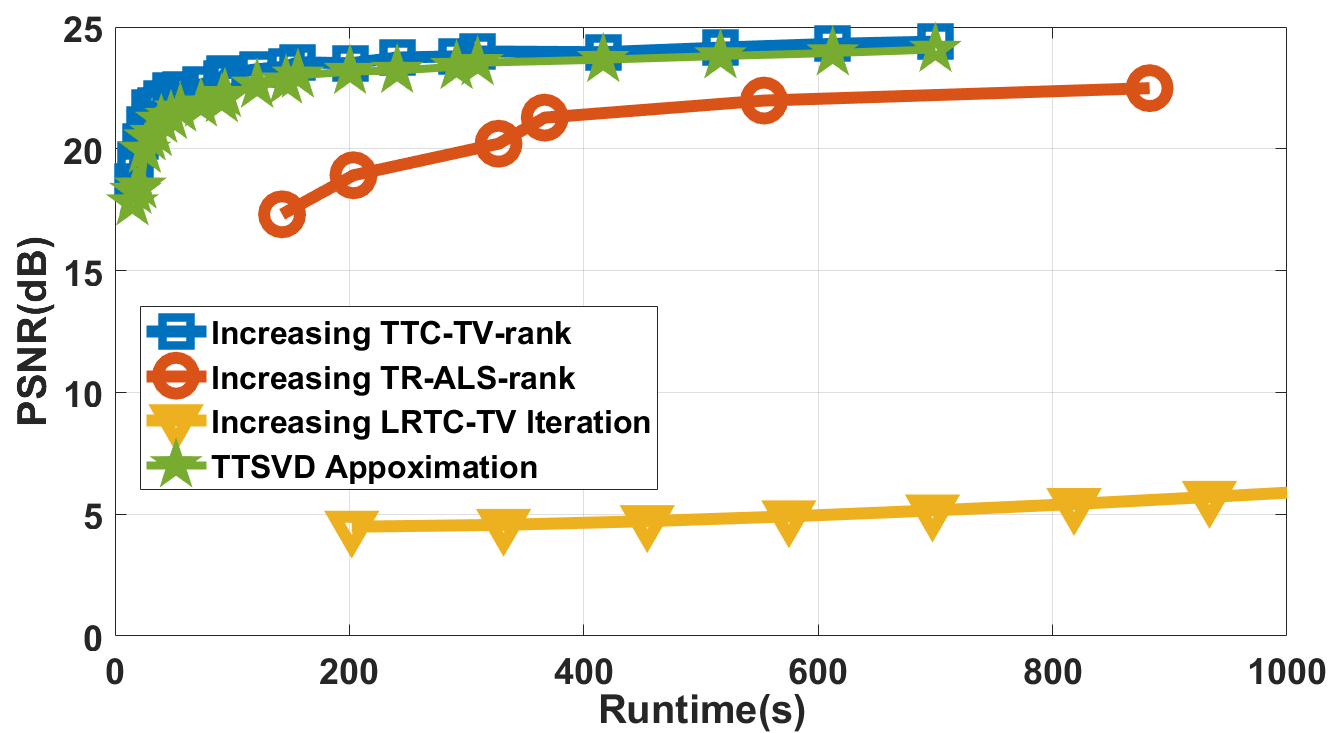}\vspace{4pt}
\end{minipage}}
\subfigure[Water Nature Fall.]{
\begin{minipage}[b]{0.45\linewidth}
\includegraphics[width=1\linewidth]{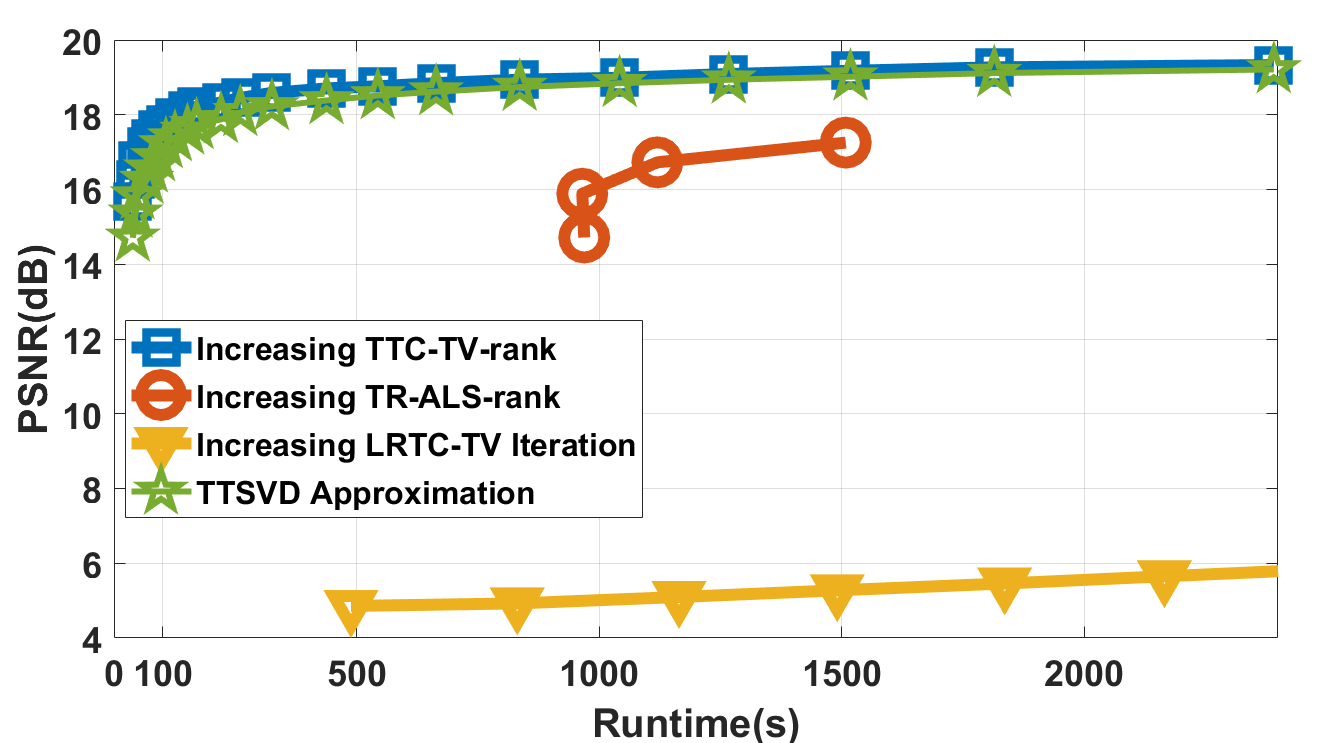}\vspace{4pt}
\end{minipage}}
\subfigure[Orion nebula.]{
\begin{minipage}[b]{0.45\linewidth}
\includegraphics[width=1\linewidth]{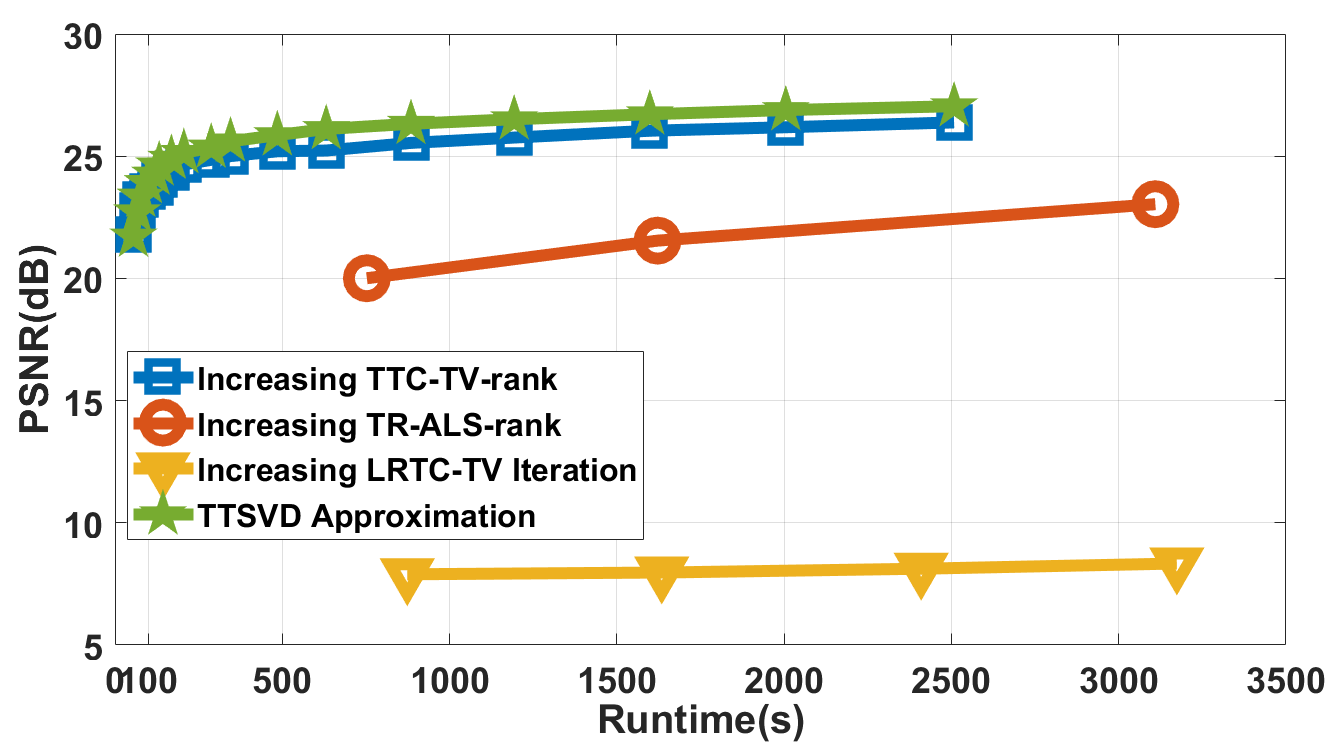}\vspace{4pt}
\end{minipage}}
\subfigure[Extended runtime Water Nature Fall.]{
\begin{minipage}[b]{0.45\linewidth}
\includegraphics[width=1\linewidth]{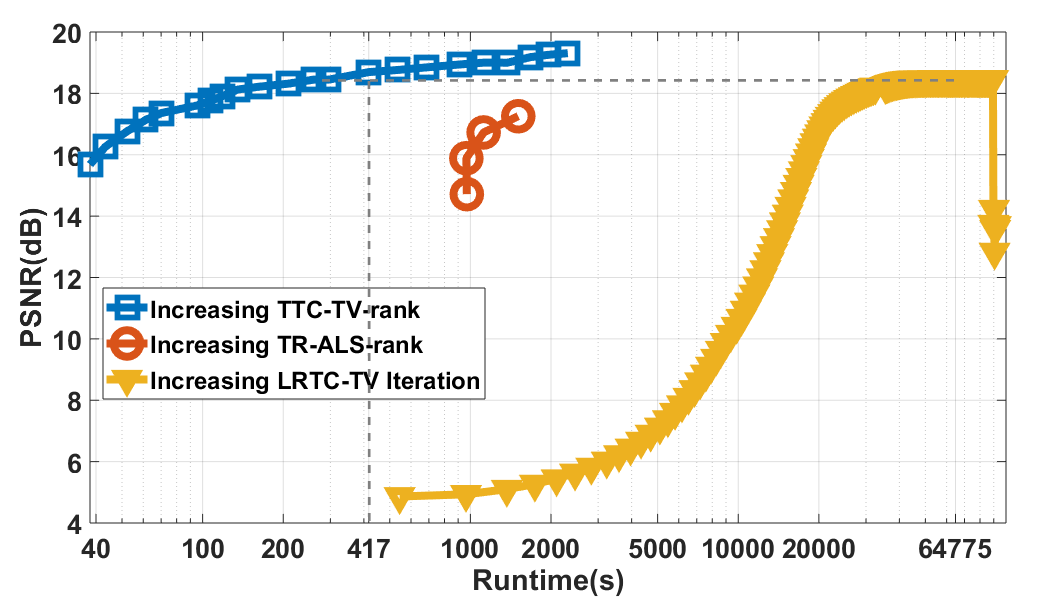}\vspace{4pt}
\end{minipage}}
\caption{Obtained PSNR of TTC-TV, TR-ALS and LRTC-TV on high-resolution benchmark images versus the total runtime.}
\label{fig:highres}
\end{figure*}
\begin{figure*}[t]
\centering
\subfigure[Original]{
\begin{minipage}[b]{0.3\linewidth}
\includegraphics[width=1\linewidth]{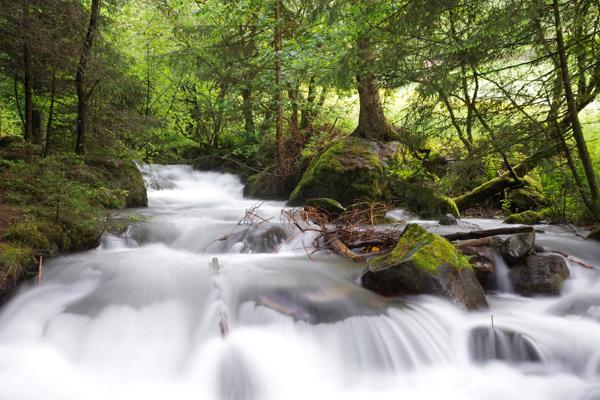}\vspace{4pt}
\end{minipage}}
\subfigure[LRTC-TV: PSNR=$18.453$]{
\begin{minipage}[b]{0.3\linewidth}
\includegraphics[width=1\linewidth]{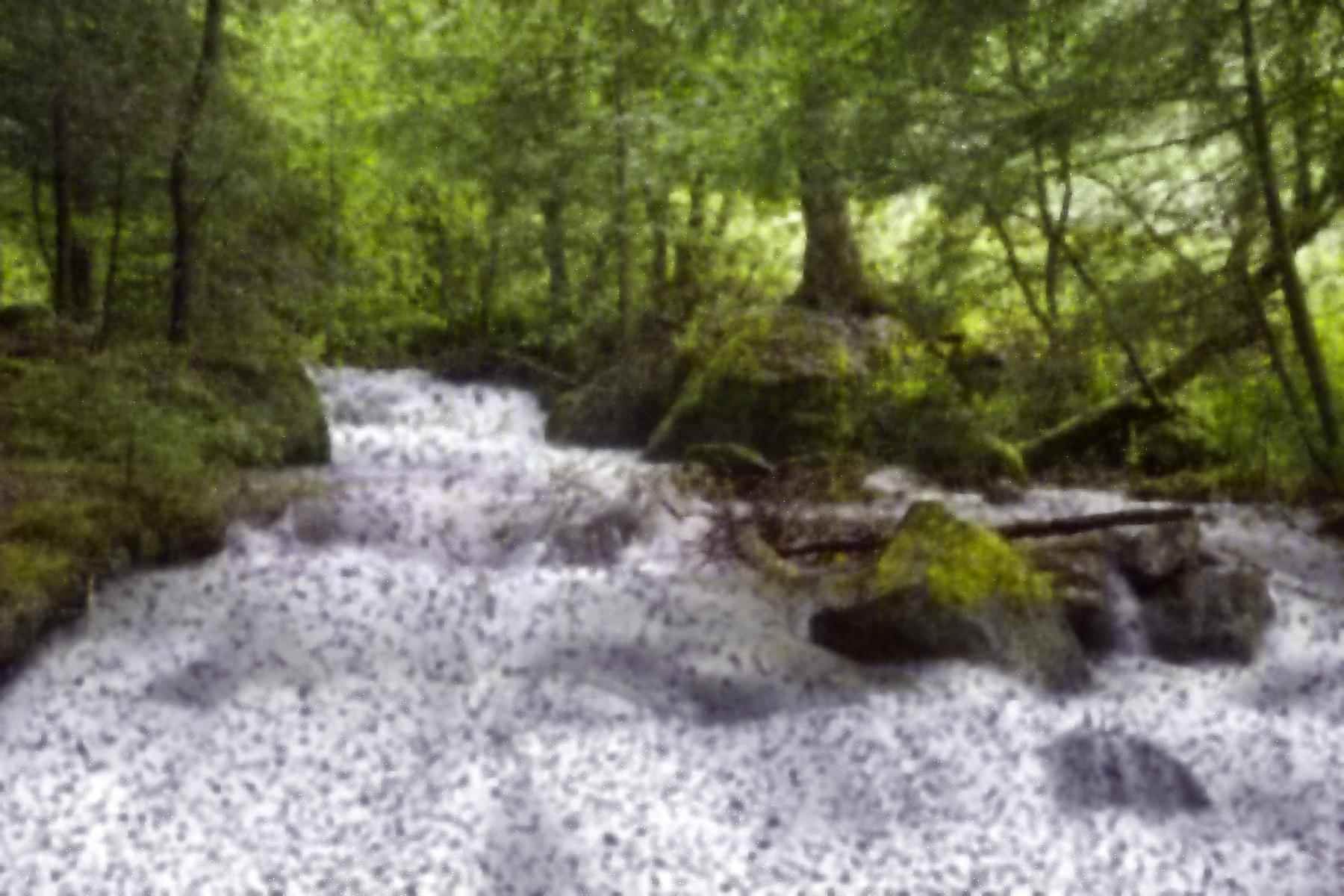}\vspace{4pt}
\end{minipage}}
\subfigure[TTC-TV: PSNR=$19.321$]{
\begin{minipage}[b]{0.3\linewidth}
\includegraphics[width=1\linewidth]{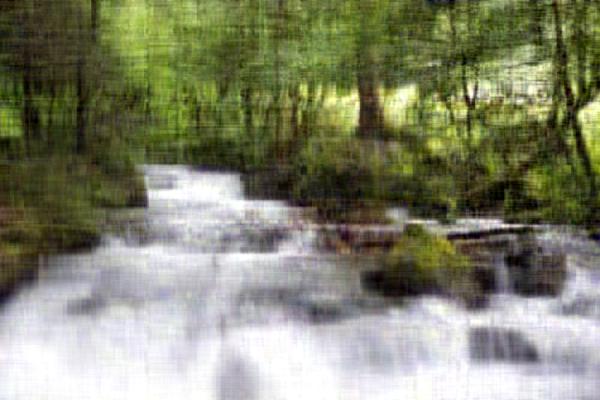}\vspace{4pt}
\end{minipage}}
\caption{Image inpaintings of \textit{Water Natural Fall} by LRTC-TV and TTC-TV.}
\label{fig:large_water_compare}
\end{figure*}

The experiments above have shown that the performance of both the TR-ALS and LRTC-TV methods are similar to the proposed TTC (TTC-TV) method. In this experiment we assess the scalability of the LRTC-TV, TR-ALS and TTC-TV methods by completing three high-resolution benchmark color images. The ground-truth of these three images is shown in the supplemental materials. The color images used in this section are \textsl{Dolphin}\footnote{http://absfreepic.com/free-photos/download/dolphin-4000$\times$3000\_21859.html}, \textsl{Water Nature Fall}\footnote{http://absfreepic.com/free\-photos/download/water-nature-fall-6000$\times$4000\_90673.html} and \textsl{Orion nebula}\footnote{http://absfreepic.com/free-photos/download/orion-nebula-in-space-6000$\times$6000\_50847.html}. The dimensions of each image and their respective factorizations used in TTC-TV and TR-ALS are listed in Table~\ref{tbl:pre2}. The LRTC-TV method uses the original dimensions of each image. The corresponding tensor trains consisted of 16, 17 and 17 cores, respectively. Only $1\%$ pixels of each image were retained. Figure~\ref{fig:highres} shows the PSNR obtained by TTC-TV, TR-ALS and LRTC-TV as a function of the total runtime for all three images. These graphs were constructed by increasing the TT-ranks for both the TTC-TV and TR-ALS methods, which resulted in better completion results at the cost of increased runtime per iteration. The TT-rank $R_{\textrm{mid}}$ of the TTC-TV method was increased from 3 up to 24, 24 and 19 for \textsl{Dolphin}, \textsl{Water Nature Fall} and \textsl{Orion nebula}, respectively. The TT-ranks could only be varied from 2 up to 7, 5 and 4 for the TR-ALS method, as the uniform TT-rank quickly lead to out-of-memory errors. This is reflected by the limited number of points in each of the TR-ALS graphs in Figure~\ref{fig:highres}. The LRTC-TV method does not use tensor trains and we allowed it to run without time restriction. 

All four figures in Figure~\ref{fig:highres} show that TR-ALS manages to achieve almost as good PSNR values as TTC-TV at about 10 times larger runtimes. The PSNR values obtained by the LRTC-TV method are unacceptable within reasonable runtimes. Figure~\ref{fig:highres}(d) illustrates that LRTC-TV needs a runtime that is about \textbf{155 times} larger (64775 versus 417 seconds) than TTC-TV to obtain the same PSNR. This is due to the use of the Tucker decomposition and a corresponding computational complexity of $O(K(\prod_{k=1}^d I_k)^3)$, which scales badly with both $d$ and $I_k$. We further note that LRTC-TV is unable to reach as high PSNR as TTC-TV in this experiment, and a sharp drop in the PSNR is witnessed after 89765 seconds as shown in Figure~\ref{fig:highres}(d). 
The completed images of LRTC-TV and TTC-TV are depicted in Figure~\ref{fig:large_water_compare}. By looking at the water flow in Figure~\ref{fig:large_water_compare}(b), it is obvious that many missing pixels remain missing after LRTC-TV completion. As the RSE curves of these experiments lead to same conclusions as above, we only show the RSE curves of \textsl{Dolphin} benchmark in
the supplemental materials. 

Moreover, we further included the curve obtained by applying TT-SVD for given TT-ranks on the original benchmark images in Figure~\ref{fig:large_water_compare}(a)(b)(c). Noted that for the given TT-ranks, the TT approximation given by TT-SVD method is quasi-optimal in the tensor-train subspace $\mathcal{S}^{(d)}_{\text{TT}}$~\cite[Corollary 2.4.]{oseledets2011tensor}. Thus including the curves enables us to evaluate the gaps between the result obtained by TTC-TV and a quasi-optimal solution in the subspace. However, we stress that the curve is only meaningful when comparing to the curve of TTC-TV, because the subspaces of TR-ALS and LRTC-TV are essentially different from that of TTC-TV.
We increased the TT-rank $R_{\textrm{mid}}$ of the TT-SVD approximation exactly the same as we did for TTC-TV, and compare the PSNRs of TTC-TV and TT-SVD approximation for same the TT-ranks. It is then seen that the result obtained by TTC-TV is nearly or even better than a quasi-optimal approximation in the tensor train subspace $\mathcal{S}^{(d)}_{\text{TT}}$.

\subsection{Color Video Completion}
\label{subsec:CVC}

\begin{figure}[t]
\centering
\includegraphics[width=0.475\textwidth]{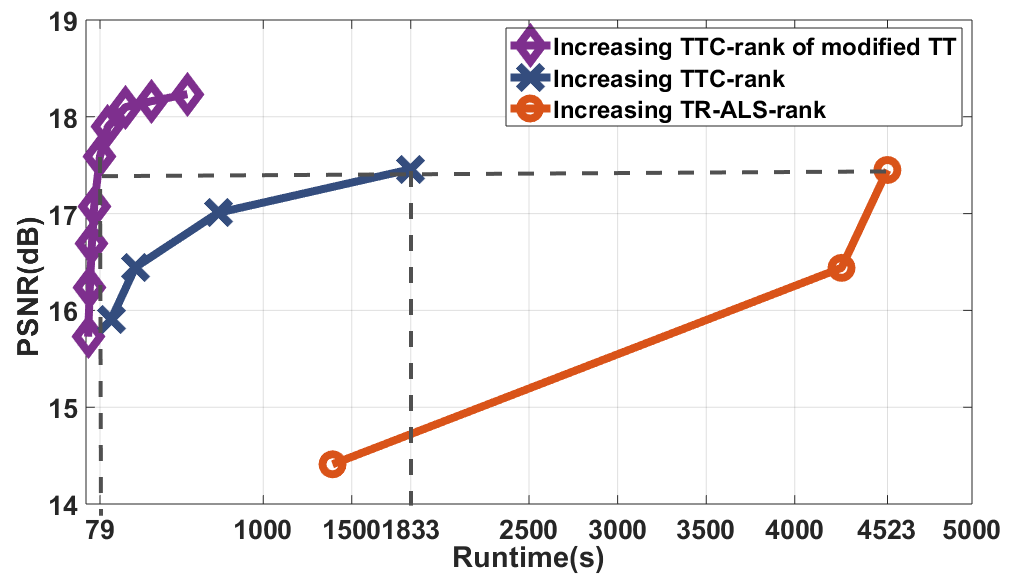}
\caption{Obtained PSNR of both TTC and TR-ALS on the video versus the total runtime for increasing TT-ranks.}
\label{fig:video}
\end{figure}

The inpainting of high-resolution images demonstrated the lack of scalability of the LRTC-TV method. This method will therefore not be considered anymore. The purpose of this experiment is to demonstrate an advantageous feature of our proposed TTC algorithm with regard to color video completion. A $144$ frames video clip was taken from the \textit{Mariano Rivera Ultimate Career Highlights} video, available on YouTube\footnote{https://www.youtube.com/watch?v=UPtDJuJMyhc}. 
The dimensions of the resulting tensor and their factorizations are given in Table~\ref{tbl:pre3}, resulting in a tensor train of eleven cores. The total number of elements of this particular tensor has the same order of magnitude as for the high-resolution images. In this experiment $10\%$ of the video is retained. Observed pixels appear consistently at the same position over all frames and color channels, which models a breakdown of $90\%$ of the available sensors in the camera. The fact that all observed pixels occur at the same position over all frames can be taken into account by solving the following optimization problem
\begin{align}
\min_{\mat{A}} \;||\mat{S}^T \mat{A} - \mat{Y}||_F^2, \textrm{ \textrm{s.t.} } \text{TT-rank}(\mat{A}) = (R_1,\ldots,R_{d}),
\label{eqn:vectoroutput}
\end{align}
where $\mat{A}$ is the original tensor reshaped into a $(360\cdot 640) \times (144\cdot 3)$ matrix. This implies that the observed values are also reshaped into the matrix $\mat{Y}$. The input matrices related to the frame and color dimensions are in this way not necessary anymore. The matrix $\mat{A}$ is then modeled by a tensor train of seven cores for which $\ten{A}^{(1)}$ has dimensions $1 \times 9 \times (144\cdot 3) \times R_2$. The $(144\cdot 3)$ dimension of $\ten{A}^{(1)}$ corresponds with the columns of the $\mat{A}$ matrix and updating this TT-core is done by rewriting~\eqref{eq:le} into a matrix equation. All other TT-cores are updated by solving~\eqref{eq:le}. Figure~\ref{fig:video} shows the PSNR as a function of the total runtime obtained by the conventional TTC algorithm, its modification~\eqref{eqn:vectoroutput}, and TR-ALS. The TT-rank for the TR-ALS method was varied from 2 up to 4, higher values for the ranks resulted in out-of-memory-errors. The best PSNR obtained with TR-ALS is $17.45$ dB and is obtained after a total runtime of 4523 seconds. The conventional TTC algorithm was run for a fixed number of 3 iterations and increasing TT-ranks from 2 up to 5 and obtains the same PSNR value about 2.5 times faster. The modified TTC algorithm for solving Equation~\eqref{eqn:vectoroutput} was run for a fixed number of 3 iterations and increasing TT-ranks from 2 up to 10. The smaller amount of TT-cores results in faster runtimes and better PSNRs. The best PSNR value obtained by TR-ALS is obtained by this particular TTC implementation in 79 seconds, about 57 times faster than TR-ALS and 23 times faster than the standard TTC formulation.

\begin{table}[t]
\scriptsize
\begin{center}
\caption{Experiment 3 dimension settings.}
\begin{tabular}{@{}rr@{}}
Original dimensions & Dimension factorization \\
\hline
$360 \times 640 \times 144 \times 3$ & $9 \times 8\times 5\times 4\times 4\times 5\times 8\times 4\times 6\times 6\times 3$ 
\end{tabular}
\label{tbl:pre3}
\end{center}
\end{table}  

\section{Conclusions}
\label{sec:C}
We have proposed an efficient tensor completion framework by assuming tensor train structures in the underlying regression model. Specifically, the multi-indices (coordinates) of the known entries act as inputs and their corresponding values act as outputs. Moreover, Total Variation regularization and Tikhonov regularization are readily realized under the tensor train framework with almost zero additional computational cost. A simple yet effective tensor train initialization method based on interpolations has also been introduced for images and videos. Extensive experiments with low percentages of known pixels have shown that the proposed algorithm not only outperforms the state-of-the-art methods in both accuracy and time cost, but also demonstrates better scalability. 

\ifCLASSOPTIONcaptionsoff
  \newpage
\fi



%
\bibliographystyle{IEEEtran}
\bibliography{TIP}

\end{document}